\documentclass[a4paper, 11pt]{amsart}
\usepackage[T1]{fontenc}
\usepackage[latin1]{inputenc}
\usepackage[english]{babel}
\usepackage{amsmath}
\usepackage{amssymb}
\usepackage{amsfonts}
\usepackage{graphicx}
\usepackage[dvips]{color}

\begin{document}
\newcommand{\normi}[1]{\|#1\|}
\newcommand{\itse}[1]{\left|\,#1\right|}
\newcommand{\its}[1]{\bigl|\,#1\bigr|}
\newcommand{\rn}{\mathbb{R}^n}
\newcommand{\na}{\mathbb{N}}
\newcommand{\re}{\mathbb{R}}
\newcommand{\R}{\mathcal{R}}
\newcommand{\Z}{\mathbb{Z}}
\newcommand{\M}{\mathcal{M}}
\newcommand{\eps}{\varepsilon}
\newcommand{\vM}{\overset{\rightarrow}{M}}
\newcommand{\ve}[1]{\overset{\rightarrow}{#1}}
\def\Xint#1{\mathchoice
 {\XXint\displaystyle\textstyle{#1}}%
 {\XXint\textstyle\scriptstyle{#1}}%
 {\XXint\scriptstyle\scriptscriptstyle{#1}}%
 {\XXint\scriptscriptstyle\scriptscriptstyle{#1}}%
 \!\int}
 \def\XXint#1#2#3{{\setbox0=\hbox{$#1{#2#3}{\int}$}
 \vcenter{\hbox{$#2#3$}}\kern-.5\wd0}}
 \def\ddashint{\Xint=}
 \def\dashint{\Xint-}
\newcommand{\m}{m}
\newcommand{\subsub}{\subset\subset}
\newcommand{\av}[3]{\underset{B(#2,#3)}{\dashint}#1(y)\,dy\,}
\newcommand{\avr}[1]{\underset{#1}{\dashint}}
\newtheorem{theorem}{Theorem}[section]
\newtheorem{lemma}[theorem]{Lemma}
\newtheorem{proposition}[theorem]{Proposition}
\newtheorem{lause}[theorem]{Theorem}
\newtheorem{definition}[theorem]{Definition}
\newtheorem{corollary}[theorem]{Corollary}
\newtheorem{question}[theorem]{Question}
\title[On the uniqueness of quasihyperbolic geodesics]
{On the uniqueness of quasihyperbolic geodesics}
\author{Hannes Luiro}
\address{Department of Mathematics and Statistics\\
University of Jyväskylä\\
P.O.Box 35 (MaD)\\
40014 University of Jyväskylä, Finland}
\email{hannes.s.luiro@jyu.fi}
\subjclass[2000]{Primary 30C65}
\keywords{Quasihyperbolic geometry, quasihyperbolic geodesic}
\maketitle
\begin{abstract}
In this work we solve a couple of well known open problems related to the quasihyperbolic metric. In the case of planar domains, our first main result states that quasihyperbolic geodesics are unique 
in simply connected domains. As the second main result, we prove that for an arbitrary plane domain $\Omega$ the geodesics between $x$ and $y$ are unique if $d_{Q}(x,y)<\pi\,$. This bound is sharp and improves the earlier bound established in \cite{Vä1}.
Concerning the $n$-dimensional case, it is shown that there exists a universal constant $c>0$ such that any quasihyperbolic ball $B_{Q}(x,r)\subset\Omega\subset\rn$ is convex if $r\leq c\,$. 
\end{abstract}

\section{Introduction}
For a domain $\Omega\subsetneq\rn$, $n\geq2$, the quasihyperbolic length of a rectifiable arc $\gamma\subset \Omega$ is defined by
\begin{equation*}
l_Q(\gamma)=\int_{\gamma}\frac{|dz|}{d(z,\partial\Omega)}\,,
\end{equation*}
where $d(z,\partial\Omega)$ is the Euclidean distance from $z$ to $\partial\Omega\,$. The \textit{quasihyperbolic metric} is defined by
\begin{equation*}
d_{Q}(x,y)=\inf l_{Q}(\gamma)\,,
\end{equation*}
where the infimum is taken over all rectifiable curves joining $x$ and $y$. The quasihyperbolic metric was introduced by F.W. Gehring and his students in the 1970's \cite{GeOs, GePa} and thereafter it has became an important tool, for example, in the theory of quasiconformal mappings. However, the precise analysis on the quasihyperbolic geometry itself has turned out to be quite challenging. In this direction, the classical results were proved by G.J. Martin and B.G. Osgood \cite{Ma}, \cite{MaOs} in the mid 1980's. In particular, they showed that quasihyperbolic geodesics are $C^1$-smooth with Lipschitz-continuous derivatives, and they also settled the quasihyperbolic geometry on the punctured plane $\re^2\setminus\{0\}$. After that, substantial progression in this particular area did not appear for about twenty years, but during the last ten years remarkable activity has been raised around the topic. For some of the further results on this field, we refer to \cite{Hä},\cite{Kl},\cite{MaVä},\cite{Vä1},\cite{Vä2}.

In this paper we concentrate on some very natural open problems, concerning the  
convexity of the quasihyperbolic balls and, especially, the uniqueness of the geodesics. The following conjectures are mentioned at least in \cite{Vä1} and \cite{Kl}:
\subsection{Uniqueness conjecture 1}\label{konj1}
There exists a universal constant $c$ such that if $x,y\in\Omega\subset\rn$, $d_{Q}(x,y)<c$, then there is only one quasihyperbolic geodesic between $x$ and $y\,$.
\subsection{Uniqueness conjecture 2}\label{konj2}
If $x,y\in\Omega\subset\re^2$ such that $d_{Q}(x,y)<\pi$, then there is only one quasihyperbolic geodesic between $x$ and $y\,$.
\subsection{Uniqueness conjecture 3}\label{konj3}
If $\Omega\subset\re^2$ is simply connected, then the geodesics are unique.
\subsection{Convexity conjecture}\label{konj4}
There exists a universal constant $c>0$ such that any quasihyperbolic ball $B_{Q}(a,r)\subset\Omega\subset\rn$ is strictly convex if  $r<c$.\\ 

The above problems were partially answered by J. Väisälä in \cite{Vä1}, where it was shown that when restricing to the planar domains, the Conjecture \ref{konj1} is valid with $c=2$ and the Convexity conjecture holds with $c=1$. As the main results of this paper, we are going to prove that 
\begin{theorem}\label{ncase}
The Convexity conjecture \ref{konj4} is valid with $c=\frac{1}{100}$ and Uniqueness conjecture 1 is valid with $c=\frac{1}{50}$.
\end{theorem}
\begin{theorem}\label{planecase}
Conjectures \ref{konj2} and \ref{konj3} are valid.
\end{theorem}
The proof of Theorem \ref{ncase} is done in Section \ref{pienetpallot}. The proof of Theorem \ref{planecase}, which turns out to be a remarkably more difficult case, is done in Sections 3-6.

We also predict that the above theorems and their proofs can be applied in different ways to develope the theory of quasihyperbolic geometry in finite and infinite dimensional spaces. However, we chose to postpone this discussion to the forthcoming studies.     

Let us introduce the outline of the proofs. For this, consider first the case $\Omega$ is convex.
In this case the above problems on uniqueness or convexity are relatively easy. For the uniqueness, for instance, it basically suffices to consider the average path of two different paths $\gamma_1,\gamma_2:[0,r]\to\Omega$ with equal length $r$ and with \textit{canonical} parametrization, that is
\begin{equation*}
|\gamma'_1(t)|=d(\gamma_1(t),\partial\Omega)\,\text{ and }\,|\gamma'_2(t)|=d(\gamma_2(t),\partial\Omega)\text{ for all }t\in[0,r]\,.
\end{equation*}
The convexity of $\Omega$ implies that $\frac{\gamma_1+\gamma_2}{2}\subset\Omega$ and especially  
the concavity of the distance function, yielding that
\begin{equation}\label{etä}
d(\frac{\gamma_1(t)+\gamma_2(t)}{2},\partial\Omega)\geq\frac{ 
d(\gamma_1(t),\partial\Omega)+d(\gamma_2(t),\partial\Omega)}{2}\,\text{ for all }t\in[0,r]\,.
\end{equation}
Combining this with the elementary fact
\begin{equation}\label{pit}
 \bigg{|}\frac{\gamma_1'(t)+\gamma_2'(t)}{2}\bigg{|}\leq \frac{|\gamma_1'(t)|+|\gamma_2'(t)|}{2}\,,
\end{equation}
guarantees that the quasihyperbolic length of the average path is less than the average of the lengths of $\gamma_1$ and $\gamma_2$. For more on the case of convex domains, we refer to \cite{Vä2} and \cite{MaVä}.

Outside the convex case, the above problems are much more difficult. Nevertheless, the proof of our first main result, Theorem \ref{ncase}, will still follow by using this idea. In short, the argument relies on a delicate estimation of how much the average path at least wins on the length element and how much it loses at most on the distance element, when the quasihyperbolic length of a path is very small. We have not tried to optimize the convexity bound ($c=\frac{1}{100}$), indeed, even the same argument with more careful error estimating would remarkably improve this bound. On the other hand, the proof reveals that the curvature of quasihyperbolic balls with small radius $r$ is uniformly comparable to $\frac{1}{r}$, just like for Euclidean balls. This fact also has a key role in the proof of Theorem \ref{planecase}. Theorem \ref{ncase} is proved in Section \ref{pienetpallot}.
   
For the desired results in the case of planar domains, Theorem \ref{planecase}, we did not succeed applying the above average-path method or its suitable modification. Instead of that, we were led to follow the framework in \cite{Vä1}, where the basis of the argument is the approximation of $\partial\Omega$ by finite sets and the use of the \textit{Voronoi diagram}. More precisely, if $\partial\Omega=\{a_1,a_2,\dots,a_N\}$, we attain for each $a\in\partial\Omega$ a \textit{Voronoi cell} $V_a$ by
\begin{equation*}
V_a=\{x\in \re^2\,:\,|x-a|= |x-\partial\Omega|\,\}.
\end{equation*}
Then each cell is a closed convex polygon (possibly unbounded) with a finite amount of \textit{corners} and \textit{edges}. The union of the cells cover the whole plane and the interiors of the cells are disjoint. In this paper the domains in $\re^2$ having finite boundary are called \textit{Voronoi domains}. 

The usefulness of the approximation by Voronoi domains is based on the fact that inside each fixed cell $V_a$ the geodesics coincide, up to a translation, with the geodesics in the punctured plane, which are known to be subarcs of \textit{logarithmic spirals}, thus images of line segments of the complex exponential mapping. This also means that the angle between $a-\gamma(t)$ and $\gamma'(t)$ is constant on $I=[t_1,t_2]$ with $\gamma_{|I}\subset int(V_a)\,$. The proof of this fact is a short calculation revealing that $e^z$ transforms the Euclidean length of a curve in $\re^2$ to the quasihyperbolic length in $\re^2\setminus\{0\}\,$. For the details, see e.g. \cite{MaOs} or \cite{Vä1}.  
 
Let us then briefly outline the main points of the proof of Theorem \ref{planecase}. Firstly, in \cite{Vä1}, the proof of the uniqueness of the geodesics for $x,y\in\Omega, d_Q(x,y)<2$, is derived directly from the obtained convexity for quasihyperbolic balls $B_Q(x,r)$ with $r<1$. In our case this method is not applicable any more, but the general idea partially resembles that. Indeed, we will show that in the case of a simply connected domain or, if $r<\pi$, the quasihyperbolic balls enjoy suitable uniform bounds for the negative curvature(see Lemma \ref{loppusjo} for the precise statement). Combining this with the fact that for small balls $B_Q(x,r)$ the curvature behaves like $\frac{1}{r}$, as it is shown in the proof of Theorem \ref{ncase}, we reach the desired uniqueness results.

The proof for the bound of the negative curvature is multistage. The first, kind of a topological part of the proof, employs the observation that two different disjoint geodesics $\gamma$ and $\tilde{\gamma}$ between $x$ and $y$ always generate a new geodesic loop strictly inside the bounded component determined by the Jordan path $\gamma\cup\tilde{\gamma}\,$. Roughly speaking, this eventually allows us to consider the existence of arbitrary narrow geodesic loops. This part of the proof is done in Section \ref{planeauxiliary}.

The most important and technically difficult part of the proof is done in Section \ref{voronoidomains}, where we show that in any Voronoi domain, if $B_Q(x,r)$ and $y\in \partial B_Q(x,r)$ are fixed, then there exists $\delta>0$ such that if $y,z\in\partial B_Q(x,r)$ and geodesics $\gamma:x\curvearrowright y$ (from $x$ to $y$), $\omega:x\curvearrowright z$ are canonically parametrized and disjoint, up to a common starting point $x$, and $|y-z|<\delta$, then 
\begin{equation}\label{introo}
\max_{t\in[0,r]}|\gamma(t)-\omega(t)|\leq e^{2r}|y-z|\,.
\end{equation}     

As the main single ingredient of the proof of (\ref{introo}) we show, 
roughly speaking, that if $\omega$ and $\gamma$ are geodesics in a fixed Voronoi domain, disjoint and close enough, then one finds a common sequence $(V_{a_k})$ of Voronoi cells for $\gamma$ and $\omega$ such that the 'angle divergence' $\beta_k-\alpha_k$ illustrated in Picture \ref{kuva:0} below, is non-decreasing with respect to $k$. In the exact proof, the signed/directed angle-functions are necessary. Recall that $\beta_k$ and $\alpha_k$ are really constants inside a fixed Voronoi cell.  
\begin{figure}[htp]
\centering
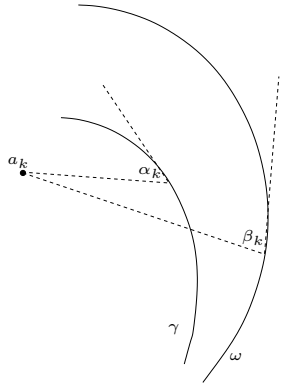
\caption{The angle divergence}
\label{kuva:0}
\end{figure}


By the above result (\ref{introo}), combined with the topological part in Section \ref{planeauxiliary}, one can prove the desired bound for the local curvature of quasihyperbolic balls on Voronoi domains simply by considering the average path of $\gamma$ and $\tilde{\gamma}$. However, to be able to achieve the desired bound for the curvature in a general domain, we need to extend these local type estimates to 'semi-global' estimates, from which the estimates in demand for general $\Omega$ can be deduced by rather straightforward approximation.
\vspace{0.5cm}

\textit{Acknowledgements.} The author is grateful to Jussi Väisälä and Riku Kl\'{e}n for their valuable comments on the manuscript.
\section{The proof of Theorem \ref{ncase}}\label{pienetpallot}

\subsection*{Canonical parametrization} 
We say that a quasihyperbolic geodesic path $\gamma:[0,r]\to\Omega$ is \textit{canonically parametrized} if $|\gamma'(t)|=d(\gamma(t),\partial\Omega)$ for all $t\in [0,r]\,$. By the definition, this also means that $r=l_Q(\gamma)$. This is the standard parametrization for the quasihyperbolic geodesics in this paper, thus if we state below that $\gamma$ is a geodesic in $\Omega$, it automatically attains the canonical parametrization for $\gamma$. For a few exceptional cases, the reader will be informed.

In this work the words 'path', 'curve' or 'arc' all refer either to the parametrizing mapping or the image of the mapping, and the meaning will be evident from the context. A rigorous treatment is given if there exists a possibility of misunderstanding. If $\gamma:[a,b]\to\Omega$ is given, then the restriction of $\gamma$ on $[c,d]\subset[a,b]$ is denoted by $\gamma_{|[c,d]}\,$. 
\begin{proposition}\label{curv}
Let $a,b\in\rn\setminus\{0\}\,$. Then
\begin{equation}\label{nojopas}
|a|+|b|-|a+b|\geq \frac{|a||b|}{2(|a|+|b|)}\bigg|\frac{a}{|a|}-\frac{b}{|b|}\bigg|^2\,.
\end{equation}
\end{proposition}
\textit{Proof.}
Observe that the right-hand-side of (\ref{nojopas}) can be written as 
\begin{equation*}
\frac{|a||b|}{2(|a|+|b|)}(2-2\frac{a}{|a|}\cdot\frac{b}{|b|})
=\frac{|a||b|-a\cdot b}{|a|+|b|}\,.
\end{equation*}

The claim follows by writing 
\begin{equation*}
|a+b|=|a|+|b|-\int_{(|a|+|b|)^2-2(|a||b|-a\cdot b)}^{(|a|+|b|)^2}\frac{1}{2\sqrt{t}}\,dt\,
\end{equation*}
to obtain that
\begin{equation*}
|a|+|b|-|a+b|\geq \frac{|a||b|-a\cdot b}{|a|+|b|}\,.
\end{equation*}
\hfill$\Box$\\

The next proposition evaluates the well known fact that for quasihyperbolic balls $B_Q(x,r)\subset\Omega$ with $r$ small enough, it holds that $d(y,\partial\Omega)\approx d(x,\partial\Omega)\,$ for all $y\in B_Q(x,r)$. We thank for Jussi Väisälä for pointing out how to improve the upper bound for $r$.
\begin{proposition}\label{quasis}
Let $n\in\na$, $\Omega\subset\rn$ a domain, $x\in\Omega$ and $r\leq\frac{1}{3}$. Then it holds that
\begin{equation*}
\frac{\sup_{y\in B_{Q}(x,r)}d(y,\partial\Omega)}{\inf_{y\in B_{Q}(x,r)}d(y,\partial\Omega)\,}=:\frac{M}{m}\leq 2\,\text{ and}\,,
\end{equation*}
\begin{equation*}
\text{if }z,y\in B_{Q}(x,r),\text{ then }\,\frac{z+y}{2}\in\Omega\text{ and } \frac{1}{2}m\leq\,d(\frac{z+y}{2},\partial\Omega)\leq 2M\,.
\end{equation*}
\end{proposition}
\textit{Proof.}
Let us denote that $\delta(a)=d(a,\partial\Omega)\,$ and recall the following classical estimates from Gehring-Palka \cite[Lemma 2.1]{GePa}: If $a,b\in \Omega$, then 
\begin{equation}\label{gehringpalka}
d_Q(a,b)\geq\log\big(\,\frac{\delta(a)}{\delta(b)}\,\big)\,\,\text{ and }\,|a-b|\leq (e^{d_Q(a,b)}-1)\delta(a)\,.
\end{equation}
Suppose then that $y,z\in B_Q(x,r)$, $r\leq\frac{1}{3}$. By the first inequality above, it follows that 
$$
\frac{\delta(z)}{\delta(y)}\leq e^{d_Q(y,z)}\leq e^{2r}\leq e^{\frac{2}{3}}<2\,,
$$
implying $\frac{M}{m}\leq 2\,$.
For the remaining claim, let us denote $w=\frac{y+z}{2}$ and use the latter inequality in (\ref{gehringpalka}) to obtain
\begin{equation*}
|w-x|\leq \max\{|y-x|,|z-x|\}\leq (e^{r}-1)\delta(x)\leq (e^{\frac{1}{3}}-1)\delta(x)<\frac{\delta(x)}{2}\,.
\end{equation*}
This implies that $|w-x|<\delta(x)$, thus $w\in\Omega\,$, and
\begin{align*}
&\delta(w)\geq \delta(x)-|w-x|>\frac{\delta(x)}{2}\geq \frac{m}{2}\,\,,\text{ and }\\
&\delta(w)\leq \delta(x)+|w-x|<\frac{3}{2}\delta(x)< 2M\,.
\end{align*}
\hfill$\Box$\\

The next proposition is a basic tool of the paper. It deals with the fact that although $d(x,\partial\Omega)$ does not need to be concave, the 'anticoncavity' is suitably bounded. 
\begin{proposition}\label{distcurv}
Let $\Omega\subset\rn$ be a domain and $x,h\in\rn$. Then 
\begin{equation*}
d(x+h,\partial\Omega)+d(x-h,\partial\Omega)\leq 2d(x,\partial\Omega)+\frac{|h|^2}{d(x,\partial\Omega)}\,.
\end{equation*}
\end{proposition}
\textit{Proof.}
Let $z,h\in\rn$. By using the elementary estimate $\sqrt{1+t}\leq 1+\frac{t}{2}$ it follows that
\begin{align*}
|z+h|&=\big(|z|^2+|h|^2+2z\cdot h\big)^{\frac{1}{2}}=|z|\big(1+\big(\frac{|h|}{|z|}\big)^2+2\frac{z\cdot h}{|z|^2}\big)^{\frac{1}{2}}\,\\
&\leq |z|\big(1+\frac{1}{2}\big(\frac{|h|}{|z|}\big)^2+\frac{z\cdot h}{|z|^2}\big)\,
=\,|z|+\frac{|h|^2}{2|z|}+\frac{z}{|z|}\cdot h\,.
\end{align*}
By applying this estimate for $h$ and $-h$, we get that 
\begin{equation*}
|z+h|+|z-h|\leq 2|z|+\frac{{|h|}^2}{|z|}\,\text{ for all }h,z\in\rn\,.
\end{equation*}
This again implies that if $a\in\partial\Omega$ and $h\in\rn$, then 
\begin{equation}\label{foinikialaiset}
d(x+h,a)+d(x-h,a)\leq 2d(x,a)+\frac{|h|^2}{d(x,a)}\,\leq 2d(x,a)+\frac{|h|^2}{d(x,\partial\Omega)}\,.
\end{equation}
Then it is easy to check that  
\begin{equation*}
d(x,\partial\Omega)=\inf_{a\in\partial\Omega}d(x,a)\,,
\end{equation*}
implies that (\ref{foinikialaiset}) is also valid if $a$ in (\ref{foinikialaiset}) is replaced with $\partial\Omega\,$.
\hfill$\Box$\\

Before proving the main theorem we state an elementary fact on the convexity:
\begin{proposition}\label{elem}
A set $G$ is strictly convex if  
\begin{equation*}
\frac{x+y}{2}\in int(G)\,\text{ for all }x,y\in \overline{G}\,.
\end{equation*}
 \end{proposition}
\vspace{2mm}
The first main results, Theorem \ref{ncase}, will be a direct corollary of the next theorem. The obtained bound for the curvature of $B_Q(x,r)$ for small $r$ will be also exploited in the planar case. Notice also that the Uniqueness conjecture 1 with constant $c=\frac{1}{50}$ follows directly from the Convexity conjecture with constant $c=\frac{1}{100}\,$.  
\begin{theorem}\label{smallballs}
Suppose that $x\in\Omega\subset\rn$, $r<\frac{1}{100}$, $y,z\in \partial B_Q(x,r)$, and 
\begin{equation*}
M:=\sup_{y\in B_{Q}(x,r)}d(y,\partial\Omega)\,.
\end{equation*}
Then 
\begin{equation*}
d_Q\big(\frac{y+z}{2}\,,\,x\,\big)\leq r-\frac{|y-z|^2}{512rM^2}\,.
\end{equation*}
\end{theorem}
\textit{Proof.}
Let $\gamma_y$ and $\gamma_z$ be geodesics from $x$ to $y$ and $x$ to $z$, respectively, and let
$$\Delta:=\sup_{t\in[0,r]}|\gamma_z(t)-\gamma_y(t)|\,.
$$ 
Consider the quasihyperbolic length of the average path 
$
\gamma=\frac{1}{2}\gamma_z+\frac{1}{2}\gamma_y\,.
$
Indeed, to prove the claim it suffices to show that
\begin{equation*}
l_Q(\gamma)\leq r-\frac{\Delta^2}{512rM^2}\,.
\end{equation*}
To show this, observe first that
\begin{align*}
l_{Q}(\gamma)-r&=\int_0^r\frac{|\gamma_y'(t)+\gamma_z'(t)|}{2d\big(\frac{\gamma_y(t)+\gamma_z(t)}{2}
,\partial\Omega\big)}\,-1\,dt\\
&=\int_0^r\frac{|\gamma_y'(t)+\gamma_z'(t)|-2d\big(\frac{\gamma_y(t)+\gamma_z(t)}{2}
,\partial\Omega\big)}{2d\big(\frac{\gamma_y(t)+\gamma_z(t)}{2}
,\partial\Omega\big)}\,dt\\
&=\int_0^r\frac{|\gamma_y'(t)+\gamma_z'(t)|-|\gamma_y'(t)|-|\gamma_z'(t)|}{2d\big(\frac{\gamma_y(t)+\gamma_z(t)}{2}
,\partial\Omega\big)}\,dt\,\\
&\,\,\,\,\,\,+\int_0^r\frac{d(\gamma_y(t),\partial\Omega)+d(\gamma_z(t),\partial\Omega)-2d\big(\frac{\gamma_y(t)+\gamma_z(t)}{2}
,\partial\Omega\big)}{2d\big(\frac{\gamma_y(t)+\gamma_z(t)}{2}
,\partial\Omega\big)}\,dt\\
&\leq 
\int_0^r\frac{-\frac{|\gamma_y'(t)||\gamma_z'(t)|}{2(|\gamma_y'(t)|+|\gamma_z'(t)|)}\bigg|\frac{\gamma_y'(t)}{|\gamma_y'(t)|}-\frac{\gamma_z'(t)}{|\gamma_z'(t)|}\bigg|^2\,}{2d\big(\frac{\gamma_y(t)+\gamma_z(t)}{2}
,\partial\Omega\big)}\,dt\,\\
&\,\,\,\,\,\,+\int_0^r\frac{|\gamma_y(t)-\gamma_z(t)|^2}{2d\big(\frac{\gamma_y(t)+\gamma_z(t)}{2}
,\partial\Omega\big)^2}\,dt\\
&=:A_1+A_2\,.
\end{align*}
Above the inequality follows directly by propositions \ref{curv} and \ref{distcurv}. For $A_2$ it follows by
Proposition \ref{quasis} (and $r<\frac{1}{3}$) that
\begin{equation*}
A_2\leq r\frac{\Delta^2}{2(m/2)^2}\,\leq 8r\frac{\Delta^2}{M^2}\,.
\end{equation*}

For the estimate of $A_1$, observe first that
\begin{align*}
&\frac{|\gamma_y'(t)||\gamma_z'(t)|}{2(|\gamma_y'(t)|+|\gamma_z'(t)|)2d\big(\frac{\gamma_y(t)+\gamma_z(t)}{2}
,\partial\Omega\big)}\\
=&\frac{d(\gamma_y(t), \partial\Omega)d(\gamma_z(t),\partial\Omega)}{2(d(\gamma_y(t), \partial\Omega)+d(\gamma_z(t))2d\big(\frac{\gamma_y(t)+\gamma_z(t)}{2}
,\partial\Omega\big)}\\
\geq & \frac{m^2}{2(2M)2(2M)}=\frac{1}{16}\bigg(\frac{m}{M}\bigg)^2\geq \frac{1}{16}\big(\frac{1}{2}\big)^2=\frac{1}{64}\,.
\end{align*}
This yields that
\begin{equation}\label{jeje}
A_1\leq -\frac{1}{64}\int_{0}^r\bigg|\frac{\gamma_y'(t)}{|\gamma_y'(t)|}-\frac{\gamma_z'(t)}{|\gamma_z'(t)|}\bigg|^2\,dt\,.
\end{equation}
For further esimates, observe that
\begin{align*}
\frac{d}{dt}(|\gamma_y(t)-\gamma_z(t)|)&=\frac{(\gamma'_y(t)-\gamma'_z(t))\cdot(\gamma_y(t)-\gamma_z(t))}{|\gamma_y(t)-\gamma_z(t)|}=:
(\gamma'_y(t)-\gamma'_z(t))\cdot\psi(t)\,\\
&=\bigg(\frac{\gamma'_y(t)}{|\gamma_y'(t)|}-\frac{\gamma'_z(t)}{|\gamma_y'(t)|}\bigg)\cdot|\gamma_y'(t)|\psi(t)\\
&=\bigg(\frac{\gamma'_y(t)}{|\gamma_y'(t)|}-\frac{\gamma'_z(t)}{|\gamma_z'(t)|}\bigg)\cdot|\gamma_y'(t)|\psi(t)\\
&\,\,\,\,\,\,\,\,+\bigg(\frac{1}{|\gamma_z'(t)|}-\frac{1}{|\gamma_y'(t)|}\bigg)|\gamma_y'(t)|\gamma_z'(t)\cdot\psi(t)\,.
\end{align*}
Then,  observe that the absolute value of the latter term above can be estimated from above by
\begin{align*}
&\bigg(\frac{1}{|\gamma_z'(t)|}-\frac{1}{|\gamma_y'(t)|}\bigg)|\gamma_y'(t)|\gamma_z'(t)|\leq \big{|}|\gamma_z'(t)|-|\gamma_y'(t)|\big{|}\\
=\,&|d(\gamma_z(t),\partial\Omega)-d(\gamma_y(t),\partial\Omega)|\,\leq |\gamma_z(t)-\gamma_y(t)|\leq \Delta\,.
\end{align*}

This implies that
\begin{align*}
 \Delta&\leq \int_{0}^{r}\bigg{|}\frac{d}{dt}(|\gamma_y(t)-\gamma_z(t)|)\bigg{|}\,dt\\
&\leq\int_{0}^r\bigg|\bigg(\frac{\gamma'_y(t)}{|\gamma_y'(t)|}-\frac{\gamma'_z(t)}{|\gamma_z'(t)}\bigg)\cdot|\gamma_y'(t)|\psi(t)\,\bigg|\,dt\,\,+r\Delta\,.
\end{align*}
Then, since $r\leq\frac{1}{2}$, $|\gamma_y'(t)|\leq M$ for all $t\in[0,r]$, and $|\psi(t)|\equiv 1$, we get from above that
\begin{equation*}
\int_{0}^r\bigg|\frac{\gamma'_y(t)}{|\gamma_y'(t)|}-\frac{\gamma'_z(t)}{|\gamma_z'(t)|}\bigg|\,dt\,\geq \frac{\Delta}{2M}\,.
\end{equation*}
Then we apply the Cauchy-Schwarz inequality to obtain that 
\begin{equation*}
\int_{0}^r\bigg|\frac{\gamma'_y(t)}{|\gamma_y'(t)|}-\frac{\gamma'_z(t)}{|\gamma_z'(t)|}\bigg|^2\,dt\,\geq \frac{1}{r}\bigg(\int_{0}^r\bigg|\frac{\gamma'_y(t)}{|\gamma_y'(t)|}-\frac{\gamma'_z(t)}{|\gamma_z'(t)|}\bigg|\,dt\,\bigg)^2\,\geq\frac{1}{r}\bigg(\frac{\Delta}{2M}\bigg)^2\,.
\end{equation*}
By combining this with (\ref{jeje}) we finally get that
\begin{equation*}
A_1\leq -\frac{1}{64}\frac{1}{r}\bigg(\frac{\Delta}{2M}\bigg)^2=-\frac{\Delta^2}{256rM^2}\,.
\end{equation*}
Summing up, it follows that
\begin{align*}
A_1+A_2&\leq -\frac{\Delta^2}{256rM^2}\,+8r\frac{\Delta^2}{M^2}\,=\frac{\Delta^2}{M^2}\big(-\frac{1}{256r}\,+
8r\big)\\
&=\frac{\Delta^2}{M^2}\big(\frac{2048r^2-1}{256r}\big) \leq -\frac{\Delta^2}{512rM^2}\,.
\end{align*}
Above the final inequality follows by $r<\frac{1}{100}\,$. This completes the proof.
\hfill$\Box$\\

Now Theorem \ref{ncase} follows from the above theorem and Proposition \ref{elem}.
   
\section{Preliminaries for the planar case}\label{planeauxiliary}
In this work we prefer to denote the underlying plane by $\re^2$ instead of $\mathbb{C}$. However,  the elementary complex analysis is utilized in the arguments. In particular, in the proof of the main Lemma \ref{major}, the use of a set-valued argument function is necessary and in this section an elementary complex integration helps us.

Notation: The imaginary unit will be denoted by $\textbf{i}$, instead of $i$, which is used for indexing. For $A\subset\re^2$, $|A|$ denotes the $2$-dimensional Lebesgue measure of $A$. For $a,b\in \re^2$, $J[a,b]$ denotes the line segment between $a$ and $b$. 

\subsection*{Angle functions}
Let $arg(z)$ denote the usual set-valued argument function on the complex plane. We define an 'signed angle'-function $ang$ for $x,y\in \mathbb{C}\setminus\{0\}$ by
\begin{equation*}
ang(x,y)=arg(x)-arg(y)\,.
\end{equation*}
The notation $Pr(a)$ is used for the principal value of $a\in\re$ modulo $2\pi$, where 
the target domain of $Pr$ is fixed to be $(-\pi,\pi]$. We will sometimes say that $x$ points \text{left} from $y$ if $0<Pr(ang(x,y)<\pi$ and, if $-\pi<Pr(ang(x,y))<0$ that $x$ points \text{right} from $y\,$.    
We will also use the notation $Pr(S)$ for sets $S\subset \re$ for which $Pr(a)=Pr(a')$ for all $a,a'\in S$. 


 

\subsection*{Curvature terminology}
We say that smooth injective path $\gamma:[0,r]\to \re^2$ is \textit{left[right] curving} on $(a,b)\subset[0,r]$ if for every $t\in[a,b)$ there exists $h_0>0$ such that 
\begin{equation*}
Pr(ang(\gamma'(t+h),\gamma'(t)))\geq 0\,\,[\,\leq\,0\,]\,,\text{ if }h\leq h_0\,.
 \end{equation*}

\vspace{2mm}
A basic tool in the proof of Theorem \ref{planecase} will be the following observation: If there exists $x,y\in \Omega$ with at least two different geodesics from $x$ to $y$, then we can find for every $\varepsilon>0$ points $\tilde{x}$ and $\tilde{y}$ with two different geodesics $\gamma,\tilde{\gamma}:[0,r]\to\Omega$ from $\tilde{x}$ to $\tilde{y}$ such that 
\begin{equation*}
|\gamma(t)-\tilde{\gamma}(t)|<\varepsilon\,\text{ for all }t\in [0,r]\,.
\end{equation*}
We will not prove exactly this result but its appropriate technical variant, settled in the following two lemmas. In order to make those lemmas more accessible, let us introduce some terminology.

\subsection*{Jordan domains}
In this paper we say that domain $\Omega\subset\re^2$ is a Jordan domain if it is bounded and there exists a Jordan curve $\gamma$ such that $\partial D=\gamma\,$. This definition relies on Jordan curve theorem, which is constantly used in the proof of the lemmas in this section. If a Jordan domain $\Omega$ is  determined by a Jordan curve $\gamma$, then we denote $\Omega=:D_{\gamma}\,$. 

\subsection*{Geodesic leaf-like loops}
Suppose that $\gamma,\tilde{\gamma}:[0,r]\to \Omega$ are geodesics from $x$ to $y$ and there exists
 $t_0\in [0,r)$ such that 
$\gamma(t)=\tilde{\gamma}(t)$ iff $t\in [0,t_0]\cup\{r\}$ and 
\begin{equation*}
\gamma_{|[0,t_0]}=\tilde{\gamma}_{|[0,t_0)}\subset \re^2\setminus \big(\bar{D}_{\gamma_{|[t_0,r]}\cup\tilde{\gamma}_{|[t_0,r]}}\big)\,.
\end{equation*}
 then we say that $\gamma\cup\tilde{\gamma}$ is a \textit{leaf-like geodesic loop}. The point $\gamma(t_0)$ (or $t_0$) will be often called \textit{the point of divergence}.

Suppose also that $\omega\subset\re^2$ is a union of a Jordan curve $\omega_1$ and 'the outer part' $\omega_2$, that is an injective continuous curve contained in $\re^2\setminus D_{\omega_1}\,$. In this case, in order to simplify notation in the forthcoming proofs, we make a convention
$D_{\omega}:=D_{\omega_1}\,$.

\vspace{2 mm}




The following lemma states that if there exists a leaf-like geodesic loop $\sigma$ between $x$ and $y$, then there exists another leaf-like geodesic loop 'strictly inside' $\sigma\,$. 
\begin{lemma}\label{topologia1}
Suppose that $\Omega\subset \re^2$ is a domain and $x,y\in\Omega$ such that there exists two quasihyperbolic geodesics $\gamma,\tilde{\gamma}:[0,r]\to \Omega$ from $x$ to $y$ such that for some $t_0\in[0,r)$ it holds that
\begin{equation}\label{coincide}
\gamma(t)=\tilde{\gamma}(t)\text{ if }0\leq t\leq t_0 \text{ and }\gamma(t)\not=\tilde{\gamma}(t) \text{ if }t_0<t<r\,,\,\text{ and }
\end{equation}
\begin{equation*}
\gamma([0,t_0))\subset \Omega\setminus \bar{D}_{\gamma\cup\tilde{\gamma}}\,.
\end{equation*} 
Then there exists (canonically parametrized) geodesics $\varrho,\tilde{\varrho}:[0,\tilde{r}]\to\Omega$ from $x$ to $\widehat{y}\in\bar{D}_{\gamma\cup\tilde{\gamma}} $ such that
\begin{equation}\label{cond1}
\varrho(t)=\tilde{\varrho}(t)=\gamma(t)=\tilde{\gamma}(t)\text{ if }t\in [0,t_0]\,,
\end{equation}
\begin{equation}\label{cond2}
\varrho(t),\,\tilde{\varrho}(t)\in \bar{D}_{\gamma\cup\tilde{\gamma}}\,,\text{ if }t\in [t_0,\tilde{r}]\,,
\end{equation}
 and there exists $t_0\leq t_{\varrho,\tilde{\varrho}}< \tilde{r}$ such that  
\begin{equation}\label{cond3}
\varrho(t)=\tilde{\varrho}(t)\text{ if }0<t\leq t_{\varrho,\tilde{\varrho}} \text{ and }\varrho(t)\not=\tilde{\varrho}(t) \text{ if }t_{\varrho,\tilde{\varrho}}<t<\tilde{r}\,,
\end{equation}
and finally, 
\begin{equation*}
|D_{\gamma\cup\tilde{\gamma}}\setminus D_{\varrho\cup\tilde{\varrho}}|>0\,.
\end{equation*}

\end{lemma}
\textit{Proof.}
Let us choose $z\in D_{\gamma\cup\tilde{\gamma}}$ arbitrarily. Suppose first that 
\begin{equation*}
d_Q(x,z)+d_Q(z,y)=d_Q(x,y)\,.
\end{equation*}
Then there exists a geodesic $\phi:[0,r]\to \Omega$ from $x$ to $y$ such that $\phi(t_z)=z$ when $t_z:=d_Q(x,z)\in(t_0,r)$. Let then 
\begin{align*}
&t^-=\inf\{t\,:\,t< t_z\,,\,\,\phi([t,t_z])\subset D_{\gamma\cup\tilde{\gamma}}\}\,,\\
&t^+=\sup\{t\,:\,t>t_z\,,\,\,\phi([t_z,t])\subset D_{\gamma\cup\tilde{\gamma}}\}\,,
\end{align*}
and define
\begin{equation*}
\varrho=
\begin{cases}
&\gamma_{|[0,t^-]}\cup\phi_{|[t^-,t^+]}\text{ if }\phi(t^-)\in\gamma\,,\\
&\tilde{\gamma}_{|[0,t^-]}\cup \phi_{|[t^-,t^+]} \text{ otherwise}\,,
       \end{cases}
\end{equation*}
and
\begin{equation*}
\tilde{\varrho}=
\begin{cases}
&\gamma_{|[0,t^+]}\text{ if }\phi(t^+)\in\gamma\,,\\
&\tilde{\gamma}_{|[0,t^+]} \text{ otherwise}\,.
       \end{cases}
\end{equation*}
Now it is easy to check that $\varrho$ and $\tilde{\varrho}$ satisfy the requirements of the claim. 

Still we have to deal with the remaining case 
\begin{equation*}
d_{Q}(x,z)+d_{Q}(z,y)>d_Q(x,y)\,.
\end{equation*}
In this case it is easy to see that there exists $\eps>0$ such that none of the geodesics from $y'\in B(y,\eps)$
to $x$ intersects $B(z,\eps)$. 
Suppose then that    
\begin{equation*}
\partial B(y,\eps)\cap \bar{D}_{\gamma\cup\tilde{\gamma}} 
\end{equation*}
is parametrized by $\nu(\lambda)$, $\lambda\in[0,1]$, such that $\nu(0)\in\gamma$ and $\nu(1)\in\tilde{\gamma}$. Notice that since $\gamma$ and $\tilde{\gamma}$ are uniformly $C^1$ and $\gamma(t)\neq\tilde{\gamma}(t)$ if $t\in(t_0,r)$, it follows that for $\eps$ small enough $\partial B(y,\eps)\cap \bar{D}_{\gamma\cup\tilde{\gamma}}$ is a connected subarc of $\partial B(y,\eps)$. 

Let then $\lambda\in(0,1)$ and suppose that $\kappa:[0,r_{\lambda}]\to \Omega$ is a (canocically parametrized) geodesic from $x$ to $\nu(\lambda)$. Then define 
\begin{equation*}
t_{\lambda}=\inf\{\,t\,:\,\kappa_{|[t,r_{\lambda}]}\subset D_{\gamma\cup\tilde{\gamma}}\}\,,
\end{equation*}
thus $\kappa(t_{\lambda})$ is the 'last point' where $\kappa$ intersects $\gamma$ or $\tilde{\gamma}$ (or both of them). Since $\nu(\lambda)$ is inside $D_{\gamma\cup\tilde{\gamma}}$, it follows that $t_0\leq t_{\lambda}<r_{\lambda}\,$. Let then 
$\omega_{\lambda}:[0,r_{\lambda}]\to\Omega$ be defined by 
\begin{equation*}
{\omega_{\lambda}}_{|[0,t_{\lambda}]}=
\begin{cases}
&\gamma_{|[0,t_{\lambda}]}\,\text{ if }\kappa(t_{\lambda})\in\gamma\,, \text{ and }\\
&\tilde{\gamma}_{|[0,t_{\lambda}]} \text{ otherwise}\,,
\end{cases}
\end{equation*}
and
\begin{equation*}
\omega_{\lambda}(t)=\kappa(t)\,\text{ if }t\in[t_{\lambda},r_{\lambda}]\,.
\end{equation*}
Especially, it follows that $\omega_{\lambda}$ is a geodesic from $x$ to $\nu(\lambda)$ contained in $\bar{D}_{\gamma\cup\tilde{\gamma}}\,$.

Let us then denote $\nu(0)=\gamma(r_0)$ and $\nu(1)=\tilde{\gamma}(r_1)$ (recall that $\nu(0)\in\gamma$, $\nu(1)\in\tilde{\gamma}$) and define for $\lambda\in(0,1)$ that
\begin{equation*}
 \sigma_{\lambda}=\gamma_{[0,r_0]}\cup\nu_{|[0,\lambda]}\cup(-\omega_{\lambda})\,,
\end{equation*}
and
\begin{equation*}
\sigma_{0}=\gamma_{|[0,r_{0}]}\cup(-\gamma_{|[0,r_{0}]})\,\text{ and }
\sigma_{1}=\gamma_{|[0,r_{0}]}\cup\nu\cup(-\tilde{\gamma}_{|[0,r_1]})\,.
\end{equation*}
Remark that, unlike usually in this paper, in the above definition of $\sigma_{\lambda}$ the direction of the parametrization of the component paths is relevant. This is due to the following definition via the complex integration over $\sigma_{\lambda}$: 
Let 
$g:[0,1]\to \re$ be defined by
\begin{equation*}
g(\lambda):=\frac{1}{2\pi}\bigg|\int_{\sigma_{\lambda}}\frac{d\rho}{\rho-z}\,\bigg|\,\,,\,\,\,\,\,
\lambda\in(0,1)\,.
\end{equation*}
It is easy to check that every $\sigma_t$ can be written as a composition of a Jordan curve and two loops, another contained in $B(y,\varepsilon)$ and another of the form $\kappa\cup(-\kappa)$, where $\kappa$ is a certain restriction of $\gamma$. Recall also our assumption that any of the geodesics $w_{\lambda}$ do not intersect $B(z,\eps)$. Combining these, by choosing $\varepsilon<<|z-y|$, it follows that the above integral over $\sigma_t$ is determined by the Jordan part, thus we deduce that $g(\lambda)\in\{0,1\}$ for every $\lambda\in[0,1]\,$. Moreover, $\gamma(0)=0$ and $\gamma(1)=1$, following by the definition of $\sigma_0$ and $\sigma_1$ above and the fact $z\in D_{\gamma\cup\tilde{\gamma}}$.

Let then
\begin{equation*}
\tilde{\lambda}:=\inf\{\lambda\in[0,1]\,:\,g(\lambda)=1\}\,
\end{equation*}
and observe that by the definition (and $\gamma(0)=0$,$\gamma(1)=1$) there exists
\begin{equation*}
\lambda_k^-\leq\tilde{\lambda}\leq \lambda_k^+\text{ such that }\lambda_k^-\to \tilde{\lambda} \text{ and }\lambda_k^+\to \tilde{\lambda}
\text{ as }k\to\infty\,,
\end{equation*}
such that 
\begin{equation*}
g(\lambda_k^-)=0\,\text{ and }g(\lambda_k^+)=1\,\,\text{ for all }k\in\na\,.
\end{equation*}
By extracting a subsequence, if needed, we may assume that $\omega_{\lambda_k^-}\to\omega$ and $\omega_{\lambda_k^+}\to\tilde{\omega}$ as $k\to\infty\,$. Since $\omega_{\lambda_k^-},\omega_{\lambda_k^+}$ are geodesics it clearly holds that $\omega$ and $\tilde{\omega}$ are geodesics from $x$ to $\nu(\tilde{\lambda})\,$. Moreover, it follows from above that 
\begin{equation*}
\int_{\gamma_{|[0,r_0]}\cup\nu_{|[0,\tilde{\lambda}]}\cup (-\omega)}\frac{d\rho}{\rho-z}\,=0
\,\,\text{ and }\,
\bigg|\int_{\gamma_{|[0,r_0]}\cup\nu_{|[0,\tilde{\lambda}]}\cup (-\tilde{\omega})}\frac{d\rho}{\rho-z}\,\bigg|\,=\,1\,.
\end{equation*}
Combining these integrals one easily deduces that 
\begin{equation*}
\bigg|\int_{\omega\cup (-\tilde{\omega})}\frac{d\rho}{\rho-z}\,\bigg|\,=\,1\,.
\end{equation*}
Especially, this means that $\omega\not\equiv\tilde{\omega}$. 

Then, let us pick $\hat{t}\in (t_0,r_{\tilde{\lambda}})$ such that $\omega(\hat{t})\neq\tilde{\omega}(\hat{t})$ (this kind of $\hat{t}$ clearly exists by $\omega\not\equiv\tilde{\omega}$) and let 
\begin{align*}
&t^+=\sup\{t>\hat{t}: \omega(t')\neq\tilde{\omega}(t') \text{ if }t'\in[\hat{t},t]\}\,,\\
&t^-=\inf\{t< \hat{t}: \omega(t')\neq\tilde{\omega}(t') \text{ if }t'\in[t,\hat{t}]\}\,.
\end{align*}
Finally, define $\varrho,\tilde{\varrho}:[0,t^+]\to\Omega$ such that 
\begin{equation*}
\begin{cases}
&\varrho(t)=\tilde{\varrho}(t)=\gamma(t)=\tilde{\gamma}(t) \text{ if } t\leq t_0,\\ 
&\varrho(t)=\tilde{\varrho}(t)=\omega(t) \text{ if } t\leq t^-\\
&\varrho(t)=\omega(t)\,\text{ and }\,\tilde{\varrho}(t)=\tilde{\omega}(t)\,\text{ if }t\in[t^-,t^+]\,. 
\end{cases}
\end{equation*}
Now it is easy to check that $\varrho$ and $\tilde{\varrho}$ satisfy the requirements of the lemma. This completes the proof. 
\hfill$\Box$\\

Next we show, by iterating the previous lemma, that any leaf-like geodesic loop $\sigma$ generates a nested sequence of leaf-like geodesic loops. 
Moreover, if we yet assume that $D_{\sigma}\cap \partial\Omega=\emptyset$, then that nested sequence of loops can be chosen to converge to a geodesic.     
\begin{lemma}\label{topologia2}
Suppose that $x,y_1\in\Omega$ with two different disjoint geodesics $\gamma_1$ and $\tilde{\gamma}_1$ from $x$ to $y_1$ such that $D_{\gamma_1\cup\tilde{\gamma}_1}\cap\partial\Omega=\emptyset$. Then there exist $\tilde{y}\not= x$, $\tilde{y}\in \bar{D}_{\gamma_1\cup\tilde{\gamma}_1}$ and a geodesic $\gamma:[0,r]\to\Omega$ from $x$ to $\tilde{y}$ so that there exists 
a sequence of points $y_k\in\Omega$, $y_k\to\tilde{y}$ as $k\to \infty$ such that for each 
$k$ there exists geodesics $\gamma_k,\tilde{\gamma}_k:[0,r_k]\to \Omega$ from $x$ to $y_k$ such that
\begin{equation}\label{req1}
\gamma\subset \bar{D}_{\gamma_k\cup\tilde{\gamma}_k}\subset \bar{D}_{\gamma_{k-1}\cup\tilde{\gamma}_{k-1}}\,\text{ for every }k\geq 2\,,
\end{equation}
\begin{equation}\label{req2}
\gamma_k(t)\to\gamma(t)\, \text{ and }\,\tilde{\gamma}_k(t)\to \gamma(t)\text{ as }k\to\infty\,,
\end{equation}
and there exists an increasing sequence of $t_k\in[0,r)$ such that
\begin{equation}\label{req3}
\gamma_k(t)=\tilde{\gamma}_k(t)=\gamma(t) \text{ if }0\leq t\leq t_k
\text{ and }\gamma_k(t)\neq\tilde{\gamma}_k(t)\text{ if }t_k<t<r_k\,.
\end{equation}
\end{lemma}
\textit{Proof.}
Let $\gamma_1$, $\tilde{\gamma}_1$ be as in the assumptions, $t_1=0$, and let $\gamma_k,\tilde{\gamma_k}$ and $t_k$, $k\in\na$, be defined inductively as follows: If $\gamma_k,\tilde{\gamma_k}$ and $t_k$ are chosen, then set $\gamma=\gamma_k$, $\tilde{\gamma}=\tilde{\gamma}_k$ and $t_k=t_0$ in Lemma \ref{topologia1} to
obtain $\varrho,\tilde{\varrho}$ and corresponding $t_{\varrho,\tilde{\varrho}}$, according the statement of the lemma. Of course there might exist many different candidates for $\varrho$ and $\tilde{\varrho}$. In this argument it is not irrelevant which pair we choose. Indeed, let us denote the set of all proper pairs $(\varrho,\tilde{\varrho})$ by $W_k$ and let 
\begin{equation*}
A_{\varrho,\tilde{\varrho}}:=|D_{\gamma_k\cup\tilde{\gamma}_k}\setminus D_{\varrho\cup\tilde{\varrho}}|\,,
\end{equation*}
for each $(\varrho,\tilde{\varrho})\in W_k\,$.
Suppose then that $(\varrho_0,\tilde{\varrho}_0)\in W_k$ is chosen such that 
\begin{equation}\label{ovelapari}
A_{\varrho_0,\tilde{\varrho}_0}\geq \frac{1}{2}\sup_{(\varrho,\tilde{\varrho})\in W_k}A_{\varrho,\tilde{\varrho}}\,.
\end{equation}
Then we set 
\begin{equation*}
\gamma_{k+1}=\varrho_0\,\text{ and }\tilde{\gamma}_{k+1}=\tilde{\varrho}_0\,,
\end{equation*}
and $t_{k+1}=t_{\varrho_0,\tilde{\varrho}_0}\,$.  
Notice that this induction is well defined since by (\ref{cond2}) and (\ref{cond3}) each $\varrho,\tilde{\varrho}$ from the statement of Lemma \ref{topologia1} really satisfy the requirements for geodesics $\gamma$, $\tilde{\gamma}$ in the statement of the lemma.

Furthermore, it follows immediately from above that 
\begin{equation}\label{req12}
\bar{D}_{\gamma_k\cup\tilde{\gamma}_k}\subset \bar{D}_{\gamma_{k-1}\cup\tilde{\gamma}_{k-1}}\,\text{ for all }k\geq 2\,.
\end{equation}
To complete the proof, it basically suffices to show that $(\gamma_k)$ and $(\tilde{\gamma}_k)$ must converge to a same geodesic. Since by (\ref{req12}) $\gamma_k\cup\tilde{\gamma}_k$ is a nested sequence of uniformly Lipschitz leaf-like geodesic loops, for the convergence of the whole sequence it suffices to prove the convergence for some subsequence\footnote{Strictly speaking, it would be enough to prove the convergence only for a subsequence}. For this, let first $\gamma_k,\tilde{\gamma}_k:[0,r_k]\to \Omega$ be the canonical parametrization of $\gamma_k$ and $\tilde{\gamma}_k$. Recall that $\gamma_k,\tilde{\gamma}_k$ are uniformly Lipschitz, $\gamma_k(r_k)=\tilde{\gamma}_k(r_k)$ and quasihyperbolic metric is uniformly bounded in $\bar{D}_{\gamma_k\cup\tilde{\gamma}_k}$. Recall also that $(t_k)$ was an increasing sequence of the points of divergence for $\gamma_k$ and $\tilde{\gamma_k}$ (defined in (\ref{req3})). Summing up, we may assume, by extracting a subsequence, that there exists Lipschitz-curves 
$\varGamma,\tilde{\varGamma}:[0,\widehat{r}]\to \Omega$ (canonically parametrized) such that 
\begin{equation*}
t_k\to \widehat{t}\,,\,\,r_k\to \widehat{r}\,,\,\,\gamma_k(r_k),\tilde{\gamma}_k(r_k)\to \varGamma(\widehat{r})=\tilde{\varGamma}(\widehat{r}) \text{ as } k\to\infty\,,
\end{equation*}
and 
\begin{equation*}
\gamma_k(t)\to \varGamma(t)\,\text{ and }\,\tilde{\gamma}_k(t)\to \tilde{\varGamma}(t) \text{ as } k\to\infty\,,\text{ if }0\leq t<\widehat{r}\,.
\end{equation*}
Since every $\gamma_k$ and $\tilde{\gamma}_k$ is geodesic, it clearly follows that both $\varGamma$ and $\tilde{\varGamma}$ must be geodesics. Notice also that by (\ref{req12})
\begin{equation}\label{ovela}
\varGamma\cup\tilde{\varGamma}\subset \bar{D}_{\gamma_k\cup\tilde{\gamma}_k} \text{ for every }k\in\na\,.
\end{equation}

For the claim we still have to show that $\varGamma(t)=\tilde{\varGamma}(t)$ for all $t\in[0,\widehat{r}]$. If $\,0\leq t\leq \widehat{t}$, this follows directly from $\gamma_k(t)=\tilde{\gamma}_k(t)$, if $0\leq t\leq t_k\to\widehat{t}$. For the case $\widehat{t}<t\leq \widehat{r}$, let 
\begin{equation*}
E=\{t\in[\widehat{t},\widehat{r}]\,:\,\varGamma(t)\neq\tilde{\varGamma}(t)\,\}\,.
\end{equation*}
If $E\neq \emptyset$ (contradicting to the claim), then  
\begin{equation*}
t^-:=\inf E\in[\widehat{t},\widehat{r}]\,.
\end{equation*}
 Since there can not exist arbitrarily small geodesic loops, we obtain the existence of $t>t^-$ such that 
$(t^-,t)\subset E$. Let us choose 
\begin{equation*}
t^+:=\sup\{\,t\,:\,[t^-,t]\subset E\}\,.
\end{equation*}
Thus, we get that $\varGamma(t)\neq\tilde{\varGamma}(t)$ if $t\in(t^-, t^+)$ and $\varGamma(t)=\tilde{\varGamma}(t)$, if $t\leq t^-$ or $t=t^+$.
This means that 
$\varGamma,\tilde{\varGamma}:[0,t^+]\to\Omega$ satisfy the assumptions of Lemma \ref{topologia1}, implying the existence of a leaf-like geodesic loop $\varGamma_2\cup\tilde{\varGamma}_2$ from $x$ to $\widehat{y}\in \bar{D}_{\varGamma\cup\tilde{\varGamma}}$ such that
\begin{equation*}
 0<c<|D_{\varGamma\cup\tilde{\varGamma}}\setminus D_{\varGamma_2\cup\tilde{\varGamma}_2}|
<|D_{\gamma_k\cup\tilde{\gamma}_k}\setminus D_{\varGamma_2\cup\tilde{\varGamma}_2}|
\end{equation*}
and $(\varGamma_2,\tilde{\varGamma}_2)\in W_k$ for all $k\in\na$. Thus, we get by the rule (\ref{ovelapari}) in the induction that 
\begin{equation*}
|D_{\gamma_k\cup\tilde{\gamma}_k}\setminus D_{\gamma_{k+1}\cup\tilde{\gamma}_{k+1}}|>
\frac{c}{2}\,,
\end{equation*}
for all $k\in\na$. This yields that $|D_{\gamma_1\cup\tilde{\gamma}_1}|=\infty$, which is the desired contradiction. Therefore, $\varGamma(t)=\tilde{\varGamma}(t)$ for all $t\in [0,\widehat{r}]$. This completes the proof.
\hfill$\Box$\\

\section{Voronoi domains}
In the study of the quasihyperbolic geometry, the approximation of an arbitrary domain by Voronoi domains (thus domains with finite boundary)  was first time used in \cite{Vä1}. 
As it was stated in the introduction, for a Voronoi domain $\Omega\subset\re^2$, $\partial\Omega=\{a_1,a_2,\dots,a_N\}\,$, the induced Voronoi cells can be defined by 
\begin{equation*}
V_{a_i}=\{x\in \re^2\,:\,|x-a_i|= |x-\partial\Omega|\,\}\,.
\end{equation*}
The alternative formulation
\begin{equation*}
V_{a_i}=\bigcap_{j\neq i}\{x\in\re^2\,:\,|x-a_i|\leq|x-a_j|\}
\end{equation*}
reveals that each Voronoi cell is a convex polygon (possibly unbounded) and the boundary of each Voronoi cell $V$ is composed of finitely many line segments (containing the endpoints) which are called \textit{edges}. The edges may intersect each others at \textit{corners}. The set of all Voronoi cells will be usually denoted by $\mathcal{V}$ and the set of all edges by $\mathcal{E}$ ( without indication to the dependence on the underlying Voronoi domain). Every Voronoi cell contains exactly one point from  $\partial\Omega$, denoted by $S_V$, and called \textit{nucleus} of $V$.

Furthermore, each edge $E$ determines exactly two $\textit{neighbours}$, that is two Voronoi cells $V$ and $\tilde{V}$ such that $E\subset V$ and $E\subset \tilde{V}$. Moreover, it is easy to check that $S_V$ and $S_{\tilde{V}}$ are located symmetry with respect to $E$.  
\vspace{2mm}

\subsection*{Geodesics on Voronoi domains}
 The \textit{local} behaviour of the geodesics in Voronoi domains is rather easy to handle. The fundamental fact (see the introduction) is that the geodesics in the interiors of the Voronoi cells can be shown to be subarcs of logarithmic spirals. The point is that for a geodesic $\gamma:[0,r]\to\Omega$, if $(a,b)\subset[0,r]$ such that 
\begin{equation*}
\gamma_{|(a,b)}\subset V\,\text{ and }\,\gamma_{|(a,b)}\cap \partial V=\{q_1,q_2,\dots,q_N\}
\end{equation*}
for some Voronoi cell $V$ of $\Omega$, then any small enough subarc of $\gamma$, which does not contain any of the points $q_i$, must be a geodesic in $\re^2\setminus S_V$, that is a subarc of a logarithmic spiral. By elementary properties of logarithmic spirals, this finally implies that the whole subarc $\gamma_{|[a,b]|}$ is a part of a logarithmic spiral w.r.t $S_V$.
 
Still notice that a fixed geodesic $\gamma:[0,r]\to\Omega$ can not have 'arbitrarily short visits' inside Voronoi cells. To be precise, denoting 
\begin{equation*}
A(t)=\{t\in[0,r]\,:\,\gamma(t)\in int(V_t)\text{ for some Voronoi cell }V_t\,\}
\end{equation*}
and for every $t\in A(t)$ that
\begin{equation*}
 t^+=\sup\{\tilde{t}\,:\,\gamma_{|[t,\tilde{t}]}\subset int(V_t)\}\,\,\text{ and }\,\,
t^-=\inf\{\tilde{t}\,:\,\gamma_{|[\tilde{t},t]}\subset int(V_t)\}
\end{equation*}
implies that
\begin{equation*}
\inf\{t^+-t^-\,:\,t\in A(t)\,\}\,>0\,.
\end{equation*}
Combining these facts and recalling that geodesics are $C^1$, it can be easily shown (see also \cite{Vä1}) that for any geodesic $\gamma$ there exists a finite composition 
\begin{equation}\label{peruskompositio}
\gamma=\bigcup_{i=1}^M\gamma_i\,,
\end{equation}
where each $\gamma_i$ is either a subarc of a logarithmic spiral with respect to a fixed $S_V$ or then a non-singular (directed) subsegment of some of the edges. The first type of subarcs will be called \textit{logarithmic parts} of $\gamma$ and the latter type of are called \textit{straight} parts of $\gamma$. If in the above composition all the elements are maximal, in the sense that $\gamma_i$ and $\gamma_{i+1}$ can lie in a same Voronoi cell only in the case they are different type of, then we say that the above composition is \textit{canonical}. 
\vspace{2mm}

The above reasoning also gives us the following result, exploited in the forthcoming proof.
\begin{proposition}\label{eirefia}
Suppose that $\gamma$ and $\tilde{\gamma}$ are geodesics in a Voronoi domain such that $\gamma(0)=\tilde{\gamma}(0)$, $\gamma'(0)=\tilde{\gamma}'(0)$, and there exist $\varepsilon>0$ and a Voronoi cell $V$ such that 
\begin{equation*}
\gamma((0,\varepsilon))\subset int(V)\,\text{ and }\,\tilde{\gamma}((0,\varepsilon))\subset int(V)\,.
\end{equation*}
Then $\gamma(t)=\tilde{\gamma}(t)$ for all $t\in (0,\varepsilon)\,$.
\end{proposition} 

\begin{figure}[htp]
\centering
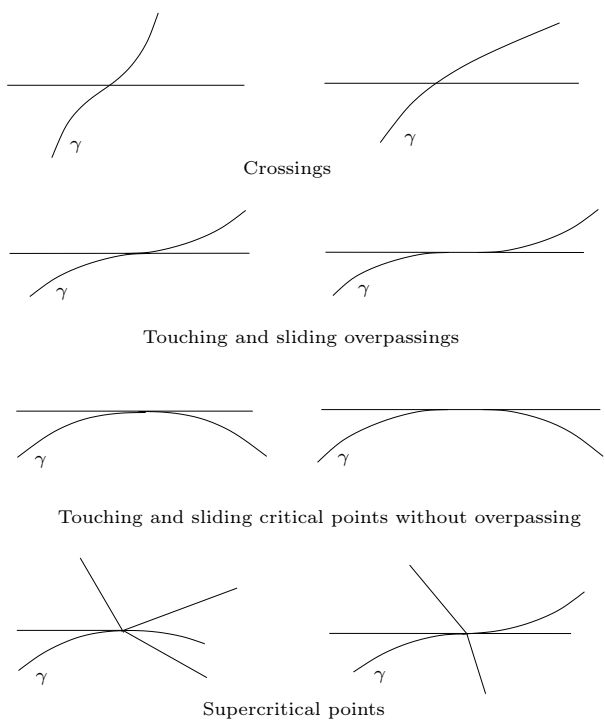
\caption{Different types of critical points}
\label{kuva:3}
\end{figure}

\subsection*{Critical points}
For a geodesic $\gamma:[0,r]\to \Omega$ in a Voronoi domain, we use to say that $a=\gamma(t)$, $t\in (0,r)$ is a \textit{critical point} of $\gamma$ if $a\in E$ for some $E\in\mathcal{E}$. Observe that if $\gamma(t)$ lies in a corner, then there exist possibly many (at least two) different edges satisfying $x\in E$. In this case, $\gamma(t)$ will be sometimes called   
\textit{supercritical} point, illustrating many technical difficulties raising up from this case.
A critical point $\gamma(t)$ will be called \textit{touching}, if $\gamma'(t)$ is parallel for some $E\in\mathcal{E}$ containing $\gamma(t)$. Furthermore, $\gamma(t)$ is called \textit{sliding} if there also exists $\varepsilon>0$ such that either $\gamma_{|[t-\varepsilon,t]}\subset E$ or $\gamma_{|[t,t+\varepsilon]}\subset E$ for some $E\in\mathcal{E}$. In the easiest case, where $\gamma(t)$ is not touching, neither lies in a corner, we say that $\gamma(t)$ is \textit{crossing}.

\section{The main lemma}\label{voronoidomains}
The following lemma contains a major part of the proof of Theorem \ref{planecase}. 
\begin{lemma}\label{major}
Let $\Omega$ be a Voronoi domain, $x,y\in \Omega$, $x\not= y$, and $\gamma:[0,r]\to \Omega$ a geodesic between $x$ and $y$. Moreover, suppose that $(\gamma_k)$ is a sequence of geodesics from $x$ to $y_k\not=y$ such that $y_k\to y$ as $k\to\infty$, $d_{Q}(x,y_k)=d_{Q}(x,y)=r$, and there exists $0\leq w_0<r$ and a non-decreasing sequence of $0\leq w_k\leq w_0$, $w_k\to w_0$ as $k\to\infty$, such that
\begin{equation}\label{helpottava}
\gamma_k(t)=\gamma(t)\text{ if and only if }t\leq w_k\,\,,\text{ and }\,|\gamma_k-\gamma|\to 0 \text{ as } k\to\infty.
\end{equation}
Then there exist $A_k\subset [0,r]$ such that $|[0,r]\setminus A_k|\to 0$ as $k\to\infty$ and
\begin{equation}\label{Crucial}
\frac{d}{dt}|\gamma(t)-\gamma_k(t)|\geq -(1+\epsilon_k)|\gamma(t)-\gamma_k(t)|\,,\text{ for all }t\in A_k\,,
\end{equation}
where $\epsilon_k\to 0$ as $k\to\infty$. Moreover, if $t\in[0,r]\setminus A_k,$ we still have, for $k$ large enough, that
\begin{equation}\label{Crucial2}
\frac{d}{dt}|\gamma(t)-\gamma_k(t)|\geq -C_{\gamma,\Omega}|\gamma(t)-\gamma_k(t)|\,.
\end{equation}
\end{lemma}

\vspace{3 mm}

Before the actual proof of the above lemma, we first prove a couple of auxiliary results.
Firstly, for expositionary reasons, we write the classical result on the regularity of quasihyperbolic geodesics (\cite{Ma}, see the above introduction) as a proposition: 
\begin{proposition}\label{apu2}
Suppose that $\Omega$ is a plane domain and $\gamma:[0,r]\to \Omega$ is canonically parametrized geodesic in $\Omega$. Then 
$$t\,\to\,\frac{\gamma'(t)}{|\gamma'(t)|}\, \text{ is }1\text{-Lipschitz on }[0,r]\,.$$ 
\end{proposition}

\begin{lemma}\label{auxiliary}
Under the conditions of the previous lemma, the following statements are valid:
\begin{enumerate}
\item ${\gamma_k}_{|[w_k,r]}$ does not intersect corners for $k$ large enough.
\vspace{1mm}

\item $\gamma_{|[w,r]}$ does not contain straight parts.
\item
\begin{equation*}\label{tokatoka}
\max_{t\in[0,r]}|\gamma_k'(t)-\gamma'(t)|\to 0\, \text{ uniformly as }k\to\infty\,.
\end{equation*} 
\item 
\begin{equation*}
\sup_{t\in (w_k,r]}\bigg{|}\big{|}Pr(ang(\gamma'(t),\gamma_k(t)-\gamma(t))\big{|}-\frac{\pi}{2}\bigg|\to 0 \text{ as } k\to\infty\,.
\end{equation*}
\item $Pr(ang(\gamma'(t),\gamma_k(t)-\gamma(t))$ is sign-preserving on $(w_k,r]$ when $k$ is large enough.
\item The last element in the canonical composition of $\gamma_k$ is logarithmic for all $k$ large enough.
\end{enumerate}
\end{lemma}
\textit{Proof.}
The proof of (1) follows directly from the facts $\gamma_k\to\gamma$ as $k\to\infty$, $\mathcal{C}$ is discrete and $\gamma_k(t)\neq\gamma(t)$ for all $t\geq w_k$. Claim (3) is also well known, the proof can be found for example from \cite[Theorem 2.8]{Vä1}. The proofs of (4) and (5) follow straightforwardly from (3). 

For the proof of (2), observe that the opposite claim implies the existence of $E\in\mathcal{E}$ and $[a,b]\subset [w,r]$, $a<b$ such that 
$
\gamma_{|[a,b]}\subset E\,.
$

Let then $V_1$ and $V_2$ be the unique Voronoi cells determined by $E$. Since $\gamma_k\to\gamma$ as $k\to\infty$ and $\gamma_k(t)\neq\gamma(t)$ for all $t\geq w_k$, it follows that for any $a<a'<b'<b$, if $k$ is large enough, then the restriction of $\gamma_k$ on $[a',b']$ is contained in $int(V_1)$ or in $int(V_2)$, thus ${\gamma_k}_{|[a',b']}$ has to coincide to the subarc of a logarithmic spiral. In particular, this clearly implies that ${\gamma_k}_{|[a',b']}$ can not be contained in an arbitrarily narrow neighbourhood of $E$. However, since $\gamma_k\to\gamma$ on $[a',b']$ and $\gamma_{|[a',b']}\subset E$, this yields the desired contradiction.

For the last claim (6), observe that the counterassumption would imply the existence of a Voronoi edge $E$ such that $\gamma(r)\in E$, $\gamma'(r)$ is parallel to $E$ and $\gamma_k(r)\in E$ for infinitely large $k$. This clearly contradicts (4) for large $k$.
\hfill$\Box$\\

Our reasoning in the proof of Lemma \ref{major} will be a sort of delicate induction with respect to common Voronoi cells/edges of $\gamma$ and $\gamma_k$. The existence of a common chain with suitable properties for our purposes is introduced in the following Lemma \ref{commoncells1}. The proof of that lemma is a bit technical, mainly due to the possible difficult type of critical points for $\gamma$.

\begin{lemma}\label{commoncells1}
Suppose that the conditions of Lemma \ref{major} are valid. Then the following holds for every fixed $k$ large enough: there exists  
\begin{align*}
&w_k=:\tilde{t}_0<\tilde{t}_1<\tilde{t}_2<\dots <\tilde{t}_{N}<\tilde{t}_{N+1}:=r\,\,,\\
&w_k=:t_0< t_1\leq t_2\leq \dots \leq t_{N}\leq t_{N+1}:=r\,\text{ and }\\
&V_0,V_1,V_2,\dots,V_{N}\in\mathcal{V}\,,
\end{align*}
such that the following hold for all $i=0,1,\dots,N$: 
\begin{enumerate}
 \item ${\gamma_k}_{|[\tilde{t}_i,\tilde{t}_{i+1}]}$ and $\gamma_{|[t_{i},t_{i+1}]}$ are contained in $V_i\,$.
\item If $i\leq N-1$, there exists an edge $E_i\subset V_i$ such that $\gamma_k(\tilde{t}_i)\in E_i$ and $\gamma(t_i)\in E_i\,$. 
\item If $V_{i+1}\neq V_i$, then $E_i$ is between $V_{i+1}$ and $V_i$. 
\item ${\gamma_k}_{|[\tilde{t}_i,\tilde{t}_{i+1}]}$ is either a logarithmic part or then a composition of a straight part and a logarithmic part (respectively).
\item If $i\geq 1$, ${\gamma}_{|[t_i,t_{i+1}]}$ is either a single point (thus $t_i=t_{i+1}$) or then a logarithmic part.
\item ${\gamma}_{|[t_0,t_1]}$ is either a logarithmic part or then a composition of a straight part and a logarithmic part (respectively).
\end{enumerate}
\end{lemma}
\textit{Proof.}
It turns out to be easier to prove the claim by contradiction. Thus, assuming the claim is not true, by extracting a subsequence, we may assume that the desired sequences do not exist for any $k\in\na$. 
Suppose then that $k$ is fixed and large enough and let 
\begin{equation}
{\gamma_k}_{|[w_k,r]}=\bigcup_{i=1}^{M}{\gamma_k}_{|[\hat{t}_i,\hat{t}_{i+1}]}=:\bigcup_{i=1}^M\omega_i\,
\end{equation}
be the canonical composition of ${\gamma_k}_{|[w_k,r]}$. Then we pick a subsequence of $(\hat{t}_i)$, thus $\hat{t}_{p_j}=:\tilde{t}_j$, $0\leq j\leq N\leq M$ according to the following inductive rule: firstly, we let $p_0=0$. Then, if $p_j$ is defined, we set $p_{j+1}$ to be either $p_j+1$ or $p_j+2$ such that the first option occurs if $\omega_{p_j+1}$ is logarithmic and the latter option if $\omega_{p_j+1}$ is straight. Since $\gamma_k$ can not have two straight parts successively (see Lemma \ref{auxiliary}), it follows that for every $0\leq j\leq N$,  ${\gamma_k}_{|[\tilde{t}_{j},\tilde{t}_{j+1}]}$ is a composition of a possible straight part and a logarithmic part, respectively. 

Since $\gamma_k$ does not intersect corners for $k$ large enough, the above choice also implies the existence of Voronoi cells $V_j$, $0\leq j\leq N$ such that 
${\gamma_k}_{|[\tilde{t}_{j},\tilde{t}_{j+1}]}$ is contained in $V_j$. Moreover, for every $1\leq j<N$ (thus $w_k<\tilde{t}_j<r$) it clearly follows that $\gamma_k(\tilde{t}_j)$ lies on edge $E_j\subset V_j$ . Moreover, if $V_j\neq V_{j-1}$, then $E_j$ is the unique edge between $V_{j-1}$ and $V_j\,$.  

Notice that above $(V_i)$, $(\tilde{t}_i)$, $N$ and $(E_i)$ all (at least \textit{a priori}) depend on $k$, unless (for simplicity) this is not indicated in the notation. However, since the amount of the elements in the canonical composition is uniformly bounded, it follows that the amount of possible outcomes for the sequences of $(V_i)$ and $(E_i)$, and $N$ above is finite. Therefore, again by extracting a subsequence, if needed, we may assume that the above operations yield fixed sequences $(V_i)$ and $(E_i)$ of fixed length $N$ for every $k\in\na$. Similarly, we may assume that related $\tilde{t}_i=:t_i^k$ (indicating the dependence on $k$) converge as $k\to\infty$ for all $1\leq i\leq N$, and thus there exists 
\begin{equation}
\lim_{k\to\infty} t^k_i\,=:\,t_i\,\text{ for all }0\leq i\leq N+1\,.
\end{equation}
Observe that since $t^k_{i+1}\geq t^k_{i}$ for every $k\in\na$ and $0\leq i\leq N$, it follows that 
$t_{i+1}\geq t_i$ for all $0\leq i\leq N\,$. Moreover, since $V_i$ and $E_i$ are closed, $\gamma_k\to\gamma$ as $k\to\infty$, and 
\begin{equation}
\gamma_k[t_i^k,t_{i+1}^k]\subset V_i\,\text{ if }0\leq i\leq N\,,\text{ and }\,\gamma_k(t_i^k)\in E_i \text{ if } 1\leq i\leq N
\end{equation}
holds for every $k\in\na$, it clearly follows that
\begin{equation}\label{fun}
\gamma[t_i,t_{i+1}]\subset V_i\,\text{ if }1\leq i\leq N\,,\text{ and }\,\gamma(t_i)\in E_i \text{ if } 0\leq i\leq N\,.
\end{equation}

For the rest of the proof, we treat separately the cases
\begin{equation}
A)\,\,w_k=w_0\text{ for all }k\in\na\,\,\,\,\,\,B)\,\,w_k<w_0 \text{ for all }k\in\na\,.
\end{equation}
By extracting a subsequence, if needed, we may assume that $A)$ or $B)$ is valid. 

In the case of $A)$, it is easy to see that the length of the first logarithmic part in the composition of $\gamma_k$ is uniformly bounded from below, and thus there exists $c>0$ such that $t^k_1>c>0$ for all $k\in\na$. This clearly implies the claim $w_k<t_1$ and moreover, 
$\gamma[w_k,t_1]=\gamma[w_0,t_1]\subset V_0$.  Summing up the above properties, we have shown that the chosen sequences $(V_i),(E_i),(t_i^k)$ and $(t_i)$ satisfy the properties (1)-(4) above, if $A)$ is valid. Furthermore, since $\gamma$ does not have straight parts on $[w_k,r]=[w_0,r]$, it follows that $\gamma[t_i,t_{i+1}]$ is logarithmic for all $0\leq i\leq N$. This implies (5) and (6) and thus we have shown the existence of the desired sequences in the case of $A)$.  
 
In the case of $B)$, observe first that $\gamma$ and $\gamma_k$ can have infinitely many different points of divergence only on the straight part of $\gamma$. Therefore, $B)$ implies that $\gamma_{|[w_k,w_0]}$ is straight for $k$ large enough and contained in a fixed edge $E$. By 
$\gamma_k(w_k)=\gamma(w_k)$ and $\gamma_k[w_k,t^1_k]\subset V_0$ for all $k\in\na$ it clearly follows that $\gamma[w_k,w_0]\subset V_0$. Since $\gamma[w_0,t_1]\subset V_0$ (see (\ref{fun}), case $i=0$), we conclude that $\gamma[w_k,t_1]\subset V_0$. Moreover, 
in this case the claim $t_1>w_k$ follows directly by $t_1=\lim_{k\to\infty}t_1^k\geq \lim_{k\to\infty}w_k=w_0>w_k\,$ for all $k\in\na$. Combining the above properties, it is shown that the chosen sequences $(V_i),(E_i),(t_i^k)$ and $(t_i)$ satisfy the required properties (1)-(4) above (if $B)$ is valid). The remaining claims (5) and (6) follow by observing that $\gamma_{|[w_k,w_0]}$ is the straight part and $\gamma_{|[w_0,t_1]}$ is logarithmic or a single point.
This completes the proof.
\hfill$\Box$\\



\subsection*{The proof of Lemma \ref{major}}  

We begin with observing that by using reflectional symmetry and a subsequence argument it suffices to prove the case where every $\gamma_k$ travels on the left-hand side of $\gamma$. That is, by using Lemma \ref{auxiliary} we may assume that
\begin{equation}\label{leftassumption}
Pr(ang(\gamma'(t),\gamma_k(t)-\gamma(t))\in (-\frac{\pi}{2}-\varepsilon,-\frac{\pi}{2}+\varepsilon)
\end{equation}
for every $k\in\na$, $t\in[0,r]\,$, and $\varepsilon>0$ can be chosen as small as wished.
Recall also that assumption (\ref{helpottava}) implies (by Lemma \ref{auxiliary}) the uniform convergence  
\begin{equation}\label{toka}
\max_{t\in[0,r]}|\gamma_k'(t)-\gamma'(t)|\to 0\, \text{ as }k\to\infty\,.
\end{equation}
The rest of the proof is divided into two parts.

\subsection*{Part 1: Non-decreasing angle divergence}
Suppose that $k$ is fixed and large enough, and thus there exist $t_i$, $t_i^k$, $0\leq i\leq N+1$, and a sequence of Voronoi cells $V_0,V_1,\dots,V_N$ and edges $E_1,E_2,\dots,E_N$ according to the statement of Lemma \ref{commoncells1} (using the notation $\tilde{t}_i=:t_i^k$). The possible  $k$-dependence is not indicated in the notation of $t_i,V_i,N$ and $E_i$. Moreover, let $J_i$ denote a line segment on $E_i$ between $\gamma(t_i)$ and $\gamma(t_i^k)$, and $\gamma(t_i^k)-\gamma(t_i)=:\vec{J}_i\,$. Recall that $S_V$ denotes the nucleus of a Voronoi cell $V$. 

Let us define for all $0\leq i\leq N$ that
\begin{equation*}
\alpha_i:=ang(\gamma'(t_{i+1}),S_{V_i}-\gamma(t_{i+1}))\,, 
\end{equation*}
\begin{align*}
\alpha_i^k:=ang(\gamma'(t^k_{i+1}),S_{V_i}-\gamma_k(t^k_{i+1})),
\end{align*}
and
\begin{equation*}
\delta_i^k:=\alpha_i-\alpha_i^k\,.
\end{equation*}
This quantity $\delta_i^k$ (or $Pr(\delta_i^k)$) is what we call the angle-divergence of $\gamma$ and $\gamma_k$. 

 \begin{figure}[htp]
\centering
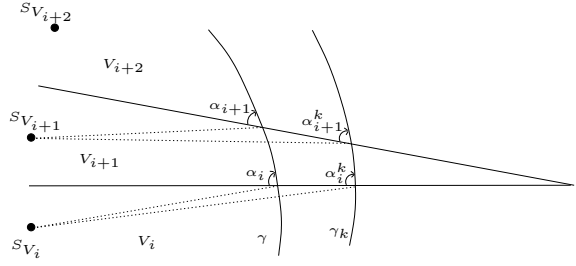
\label{kuva:15}
\caption{Angle divergence}
\end{figure}


Next we are going to show by induction that  
\begin{equation}\label{claim}
Pr(\delta^k_i)\in(-\delta,0),
\end{equation}
for all $0\leq i\leq N\,$ and $k\geq k_0$, where $\delta>0$ can be chosen as small as desired by choosing $k_0$ large enough. Simultaneously, it will be shown that 
\begin{equation}\label{gammacommon2}
{\gamma_k}_{|[t^k_i,t^k_{i+1}]}\text{ is the logarithmic part }
\end{equation}
for every $1\leq i\leq N\,$.


Notice that (by Lemma \ref{commoncells1}) it holds that $\gamma_{|[t_i,t_{i+1}]}$ is logarithmic (or a single point) w.r.t. $S_{V_i}$ for every $1\leq i\leq N$. Therefore, $\alpha_i$, $i\geq 1$ could be equivalently defined by 
\begin{equation}\label{gammacommon3}
\alpha_i=ang(\gamma'(t),S_{V_i}-\gamma(t))\,, \text{ if } t\in[t_{i},t_{i+1}]\,.\\
\end{equation}
For $\gamma_k$, we know at this point that if (\ref{gammacommon2}) holds for $\gamma_k$, then
\begin{equation}\label{keskinen}
\alpha^k_i=ang(\gamma_k'(t),S_{V_i}-\gamma_k(t))\,, \text{ if } t\in[t^k_{i},t^k_{i+1}]\,.
\end{equation}
Therefore, if (\ref{gammacommon2}) holds for $\gamma_k$, then the angle divergence can be also illustrated as in the picture 4 below.

\begin{figure}[htp]\label{kuva:15}
\centering
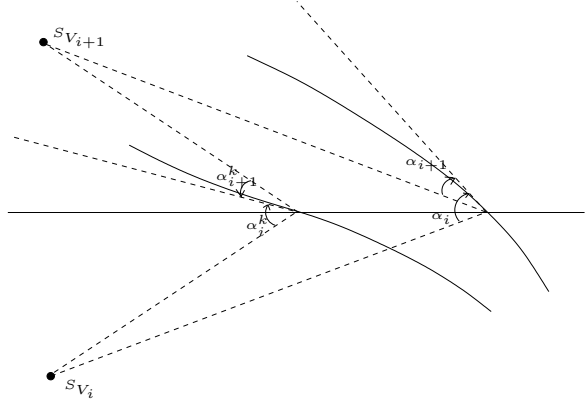
\label{abraka}
\caption{Angle divergence II}
\end{figure}

Let us then begin the inductive proof of (\ref{claim}) and (\ref{gammacommon2}).
For this, suppose that $0\leq i\leq N-1$ and 
$Pr(\delta^k_i)\in(-\delta',0)$. Consider first (\ref{gammacommon2}) for $i+1$: Suppose that (\ref{gammacommon2}) does not hold for $i+1$. This means (by Lemma \ref{commoncells1}) that $\gamma_k$ is straight on $[t_{i+1}^k,t_{i+1}^k+h]$ for some $h>0$. Applying $\gamma'_k\to\gamma'$ as $k\to\infty$, this can happen only if $\gamma'(t_{i+1})$ is also parallel to $E_{i+1}$ and points in the same direction as $\gamma_k'(t_{i+1}^k)$. 
By suitable rotation and translation (which both preserve the directed angles), it then suffices to deal with the case 
\begin{equation}\label{sup1}
V_{i+1}\subset \{\,z\,:\,Re(z)\geq 0\,\}\,,\,\,E_{i+1}\subset\{\,z\,:\,Re(z)=0\,\}\,\text{ and }\gamma(t_{i+1})=0\,.
\end{equation}
Observe that these assumptions imply that
\begin{equation}\label{sup3}
S_{V_i}\in \{\,z\,:\,Re(z)<0\,\}\,\text{ and }\,\gamma_k(t^k_i)\in \{\,z\,:\,Re(z)=0\,,\,\,Im(z)> 0\,\}\,,
\end{equation}
where the latter claim follows easily by the left-hand side assumption (\ref{leftassumption}).
The desired contradiction follows then by combining the above with  
$Pr(\delta^k_i)\in (-\delta',0)\,$. Indeed, it follows that 
\begin{align*}
Pr(\delta_i^k)&=Pr[ang(\gamma'(t_{i+1}),S_{V_i}-\gamma(t_{i+1}))-
ang(\gamma'(t^k_{i+1}),S_{V_i}-\gamma_k(t^k_{i+1}))\,]\\
&=Pr[ang(\gamma'(t_{i+1}),S_{V_i})-\gamma(t_{i+1})-
ang(\gamma'(t_{i+1}),S_{V_i}-\gamma_k(t^k_{i+1}))\,]\\
&=Pr[ang(S_{V_i}-\gamma_k(t^k_{i+1}),S_{V_i})]\,>0\,. 
\end{align*}
This contradicts with $Pr(\delta^k_i)\in (-\delta',0)\,$. Thus we have shown that $Pr(\delta^k_i)\in(-\delta',0)$, $0\leq i\leq N-1$, implies (\ref{gammacommon2}) for $i+1$.


Because (\ref{gammacommon2}) holds for $i+1$, we get that
\begin{align}
\notag\alpha_{i+1}^k&=ang(\gamma'(t^k_{i+2}),S_{V_{i+1}}-\gamma(t^k_{i+2})\\
&=ang(\gamma'(t^k_{i+1}),S_{V_{i+1}}-\gamma(t^k_{i+1}))\,.\label{keskinen2}
\end{align}
For $\gamma$ the corresponding property is already known by (\ref{gammacommon3}). 
Consider then the other part of the induction, indeed, let us show that $Pr(\delta_{i+1}^k)\in (-\delta,0)$. Let us first treat the (curious) case $V_{i+1}=V_i$. Indeed, this case was not excluded in Lemma \ref{commoncells1}. However, from the proof of that lemma one could easily see that $V_{i+1}=V_i$ can occur only if $\gamma_k$ is not logarithmic on $[t_{i+1}^k,t_{i+2}^k]\,$, and thus it would follow from above that actually $V_{i+1}\neq V_i$. Alternatively, one can also handle this case by observing that since $\gamma_k$ is logarithmic on   
$V_{i+1}$ and $S_{V_i}=S_{V_{i+1}}$, it follows directly from the definition (and (\ref{keskinen2})) that $\alpha_{i+1}^k=\alpha_{i}^k$ and $\alpha_{i+1}=\alpha_i$, and thus $Pr(\delta_{i+1}^k)=Pr(\delta_i^k)>0\,$.

For the more important case $V_{i+1}\neq V_i$, we write using the first part of the induction that 
\begin{align*}
&\alpha_{i+1}-\alpha_{i+1}^k\\=&
\,\alpha_{i}-\alpha_{i}^k+(\alpha_{i+1}-\alpha_i)-(\alpha_{i+1}^k-\alpha_i^k)\\
=& \,\alpha_{i}-\alpha_{i}^k+ 
\, ang(\gamma'(t_{i+1}),S_{V_{i+1}}-\gamma(t_{i+1}))-ang(\gamma'(t_{i+1}),S_{V_i}-\gamma(t_{i+1}))\\
&- ang(\gamma_k'(t_{i+1}^k),S_{V_{i+1}}-\gamma_k(t^k_{i+1}))+ang(\gamma_k'(t_{i+1}^k),S_{V_i}-\gamma_k(t_{i+1}^k))\\
=&\,\alpha_i-\alpha_i^k+ arg(S_{V_i}-\gamma(t_{i+1}))-arg(S_{V_{i+1}}-\gamma(t_{i+1}))
\\&-arg(S_{V_i}-\gamma_k(t_{i+1}^k))+arg(S_{V_{i+1}}-\gamma_k(t_{i+1}^k))\\
=&\,\alpha_i-\alpha_i^k+ang(S_{V_i}-\gamma(t_{i+1}),S_{V_i}-\gamma_k(t_{i+1}^k))
\\&-ang(S_{V_{i+1}}-\gamma(t_{i+1}),S_{V_{i+1}}-\gamma_k(t_{i+1}^k))\,\\
=&\,\alpha_i-\alpha_i^k+2\,ang(S_{V_i}-\gamma(t_{i+1}),S_{V_i}-\gamma_k(t_{i+1}^k))\,.
\end{align*}
The last equality is valid, because $S_{V_i}$ and $S_{V_{i+1}}$ are located symmetrically with respect to the common edge $E_{i+1}$ (as it was mentioned in the previous section). This is the first key ingredient of the proof.

Moreover, it follows that
\begin{equation*}
 Pr(ang(S_{V_i}-\gamma(t_{i+1}),S_{V_i}-\gamma_k(t_{i+1}^k)))\in(-\frac{\pi}{4},0),
\end{equation*}
when $k$ is large enough. To see this, observe first that the lower bound $\frac{\pi}{4}$ is automatically true. We also remark that it is required here for purely technical reasons, ensuring that the arguments of $Pr$ in the induction stay in the correct range. The upper bound is the second key ingredient of the proof. To show it, notice that again by applying a suitable translation and rotation we may assume that $E_{i+1}$ is contained in the imaginary-axis, $re(S_{V_i})<0$, and $\gamma(t_{i+1})=0$. As before, (\ref{leftassumption}) implies that $im(\gamma_k(t_{i+1}^k))$ must be positive. Thus we have reduced the claim to 
 \begin{align*}
 &Pr(ang((b_1+\textbf{i}b_2),(b_1+\textbf{i}(b_2-\lambda)))\\
=&Pr(arg(b_1+\textbf{i}b_2)-arg(b_1+\textbf{i}(b_2-\lambda)))\in(-\pi,0),
\end{align*}  
if $\lambda>0$, $b_1<0$ and $b_2\in\re\,$. This claim is clearly true. 

From above we conclude that to complete the proof it suffices to show that 
\begin{equation}\label{ekavaihe}
Pr(\delta^k_0)\in (-\delta,0)\,.
\end{equation}
To become convinced on the validity of this claim, see Picture \ref{firstcell2} below, which makes the claim seem rather obvious. However, we chose to express a rigorous proof at this point. 
First of all, by the left-hand side assumption (\ref{leftassumption}), it directly follows that 
\begin{equation}\label{leftassumption:seuraus}
Pr[ang(\gamma'(w_k),\gamma'_k(w_k))]\in [-\delta,0]\,.
\end{equation}

\begin{figure}[htp]
\centering
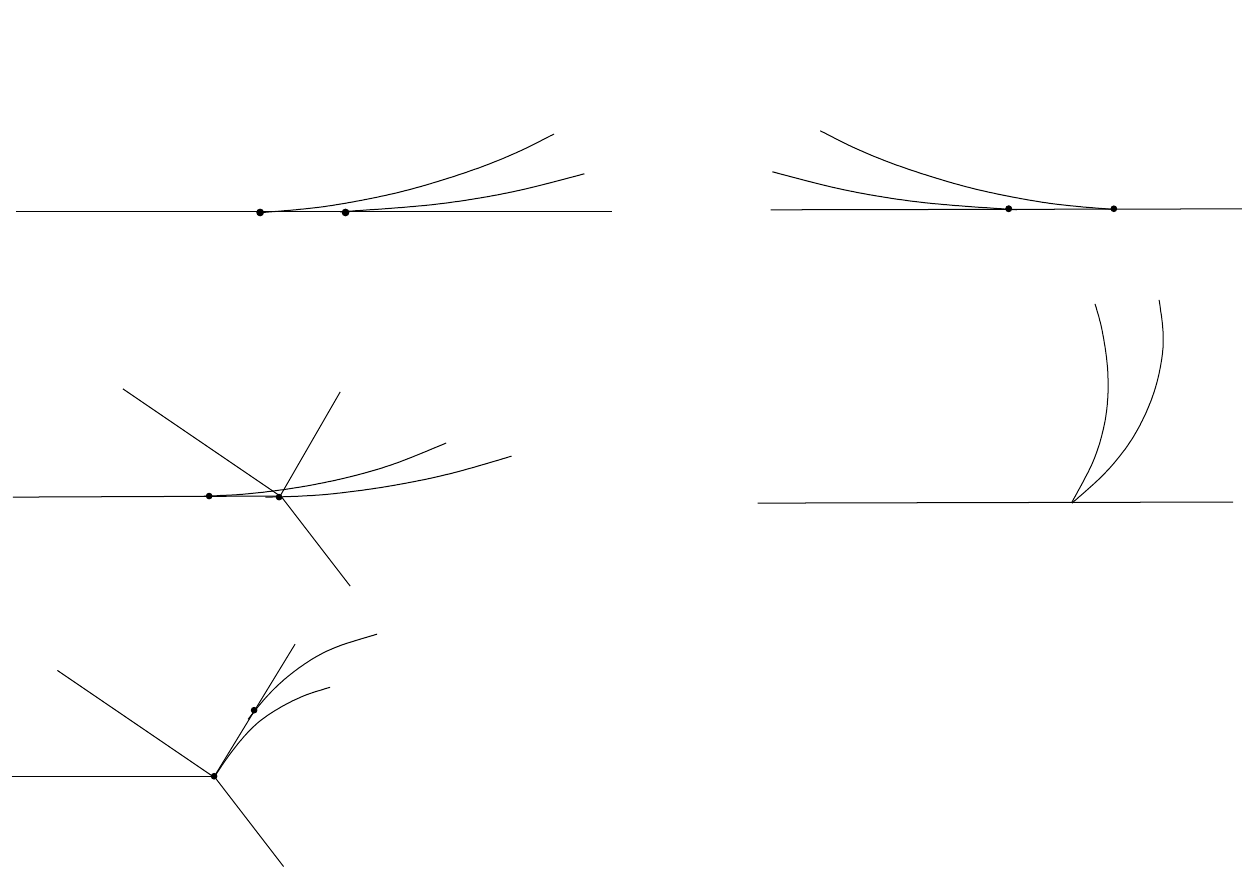
\caption{}
\label{firstcell2}
\end{figure}

\subsubsection*{Proof of (\ref{ekavaihe}), case 1:} 
Assume that there does not exist $E\in\mathcal{E}$ such that $\gamma(w_0)\in E$ with $\gamma'(w_0)$ parallel to $E$. Especially, this holds if $\gamma(w_0)\in int(V_0)\,$. In this case it follows that $w_k=w_0=0$ and also that $\gamma_{|[w_k,t_1]}$ and ${\gamma_k}_{|[w_k,t^k_1]}$ are logarithmic.
Therefore, we have
\begin{align*}
Pr(\delta^k_0)&=Pr[ang(\gamma'(t_{1}),S_{V_0}-\gamma(t_{1}))-ang(\gamma'(t^k_{1}),S_{V_0}-
\gamma(t^k_{1}))]\\
&=Pr[ang(\gamma'(w_k),S_{V_0}-\gamma(w_k))-ang(\gamma_k'(w_k),S_{V_0}-\gamma_k(w_k))]\\
&=Pr[ang(\gamma'(w_k),\gamma_k'(w_k)]\in (-\delta,0)\,.
\end{align*}

\subsubsection*{Proof of (\ref{ekavaihe}), case 2:}
This case is illustrated in Picture 6 below. Indeed, suppose that $w_k=w_0$ if $k\geq k_0$ and assume that $\gamma(w_0)=\gamma_k(w_k)\in E\in\mathcal{E}$ and $\gamma'(w_0)$ is parallel to $E$. As before, it again holds that 
$\gamma_{|[w_0,t_1]}$ is logarithmic (recall $t_1>w_k=w_0$ from Lemma \ref{commoncells1}). 
If ${\gamma_k}_{|[w_k,t^k_1]}={\gamma_k}_{|[w_0,t^k_1]}$ is also logarithmic (this can actually happen only if $w_0=w_k=0$), (\ref{ekavaihe}) follows by exactly the same argument as above. In the opposite case, it is easy to see that ${\gamma_k}_{|[w_k,t^k_1]}$ is a composition of a short straight part on $E$ and the logarithmic part, that is 
\begin{equation*}
{\gamma_k}_{|[w_k,t^k_1]}={\gamma_k}_{|[w_0,t^k_1]}={\gamma_k}_{|[w_0,w_0+\varepsilon_k]}\cup{\gamma_k}_{|[w_0+\varepsilon_k,t^k_1]}\,,
\end{equation*}
where ${\gamma_k}_{|[w_0,w_0+\varepsilon_k]}\subset E$ and ${\gamma_k}_{|[w_0+\varepsilon_k,t^k_1]}$ is a logarithmic subarc in $V_0$. Then, by a suitable rotation and translation, we may assume that $E$ lies on the real-axis, $V_0\subset\{z:im(z)\geq 0\}$ and $0=\gamma(w_0)=\gamma_k(w_0)\in E$. By the left-hand side assumption it follows that $\gamma'(w_0)$ and ${\gamma_k}_{|(w_0,w_0+\varepsilon_k]}$ are contained
in the negative real-axis. Observe also that $\gamma_k'(w_0+\varepsilon_k)$ and $\gamma'(w_0)$ point to the same direction. As in the previous case, we then have
\begin{align*}
&Pr(\delta^k_0)\\
=&Pr[ang(\gamma'(t_{1}),S_{V_0}-\gamma(t_{1}))-ang(\gamma'(t^k_{1}),S_{V_0}-
\gamma(t^k_{1}))]\\
=&Pr[ang(\gamma'(w_0),S_{V_0})-ang(\gamma_k'(w_0+\varepsilon_k),S_{V_0}-\gamma_k(w_0+
\varepsilon_k))\,]\\
=&Pr[ang(\gamma'(w_0),S_{V_0})-ang(\gamma'(w_0),S_{V_0}-\gamma_k(w_0+
\varepsilon_k))\\
=&Pr[arg(S_{V_0}-\gamma_k(w_0+\varepsilon_k))- arg(S_{V_0})]\,.
\end{align*}
Finally, this quantity is negative since $S_{V_0}$ lies strictly inside the upper half-plane and $\gamma_k(w_0+\varepsilon_k)$ on the negative real-axis. 

\begin{figure}[htp]
\centering
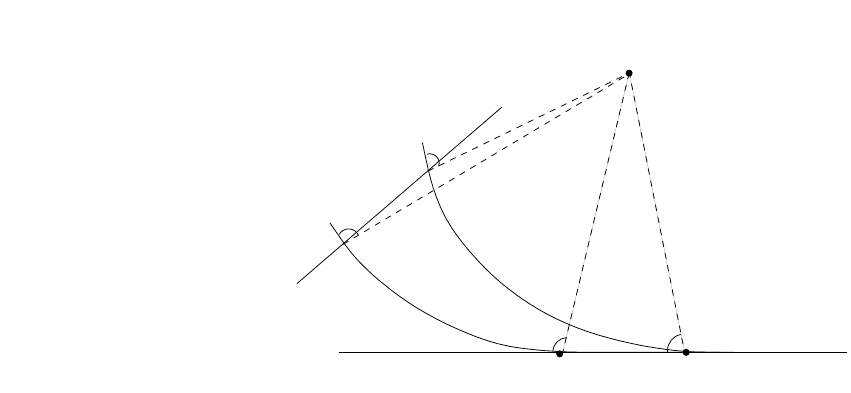
\label{kuva:7}
\caption{Case 2}
\end{figure}

\subsubsection*{Proof of (\ref{ekavaihe}), case 3:}
Assume then that $0<w_k<w_0$ for all $k\in\na$. Recalling that this means that $\gamma(w_k)$ is a point of divergence for $\gamma$ and $\gamma_k$, it is easy to see that $0<w_k<w_0$ for all $k\in\na$ combined with $w_k\to w_0$ as $k\to\infty$ implies that $\gamma$ must have a straight part 'before' the limit of the divergence points. That is, there exists $E\in\mathcal{E}$ and $\tilde{w}<w_0$ such that $\gamma_{|[\tilde{w},w_0]}\subset E\in\mathcal{E}\,$. On the other hand, now ${\gamma_k}_{|[w_k,t^k_1]}$ does not have a straight part, and thus it is a subarc of a logarithmic spiral in $V_0$. 
Therefore, 
by a similar reasoning as before, we get that
\begin{align*}
Pr[\delta^k_0]&=Pr[ang(\gamma'(w_0),S_{V_0}-\gamma(w_0))-ang(\gamma_k'(w_k),S_{V_0}-\gamma_k(w_k))]
\\
&=Pr[ang(\gamma'(w_0),S_{V_0}-\gamma(w_0))-ang(\gamma'(w_0),S_{V_0}-\gamma_k(w_k))]\\
&=Pr[arg(S_{V_0}-\gamma_k(w_k))-arg(S_{V_0}-\gamma(w_0))]\,.
\end{align*}
Again, by a suitable translation and rotation, we may assume that $V_0$ and $S_{V_0}$ lie in the upper half-plane, $\gamma_k(w_k)=0=\gamma(w_k)$ and $E$ is contained in the real axis. Recalling the left-hand side assumption, we get that ${\gamma}_{|[w_k,w_0]}$ must be contained in the positive real axis. Therefore, it follows that 
\begin{align*}
&Pr[arg(S_{V_0}-\gamma_k(w_k))-arg(S_{V_0}-\gamma(w_0))]\\
=&Pr[arg(S_{V_0})-arg(S_{V_0}-\gamma(w_0))]\,<0\,. 
\end{align*} 
This proves the case $0<w_k<w_0$.
\smallskip

Summing up the cases 1-3, we have shown (\ref{ekavaihe}), and thus the proof of (\ref{claim}) and (\ref{gammacommon2}) is complete. 


 

\subsection*{Part 2: Conclusion}
In this part we are going to show the proper claim of this Lemma, (\ref{Crucial}) and (\ref{Crucial2}). For that, let us pick from $(V_0,V_1,\dots,V_N)$ those $V_i$ for which $t_i<t_{i+1}$. 
Thus, we get a subsequence $(V_{p_1}, V_{p_2},\dots, V_{p_l})\,$, $0\leq p_1< p_2<\dots< p_{l}\leq N$ such that $t_{p_j}<t_{p_j+1}\,$ for all $1\leq j\leq l\,$. Observe that if $\gamma(t_{p_j+1})$ lies in a corner, then it might happen that $p_{j+1}>p_j+1$. It also holds that  $t_0=t^k_0=w_k<t_1$ (by Lemma \ref{commoncells1}), implying that $p_1=0$. 

Let then
\begin{align*}
a^k_{j}:=\max\{t_{p_{j}},t^k_{p_{j}}\}\,\,\text{ and }
\,\,b^k_{j}:=\min\{t_{p_j+1},t^k_{p_j+1}\}.
\end{align*}
It is easy to check that if $k$ is big enough, then 
\begin{equation*}
a^k_j<b^k_j\,\,\text{ for all }j=1,2\dots, l\,.
\end{equation*}
Moreover, if this is the case, then it follows by the definition that 
\begin{equation*}
\gamma_k(t),\gamma(t)\in V_{p_j}\,\text{ if }t\in [a^k_j,b^k_j]\,.
\end{equation*}
Yet notice that $t_0=t^k_0=w_k<t_1$ and $p_1=0$ imply that $a^k_1=w_k\,$.

\subsubsection*{The proof of (\ref{Crucial})}
Let $1\leq j\leq l$ and denote the nucleus of $V_{p_j}$ by $S$. We write that
\begin{align*}
&ang(\gamma_k'(t),\gamma'(t))=arg(\gamma_k'(t))-arg(\gamma'(t))\\=&
arg(\gamma_k'(t))-arg(S-\gamma_k(t))+arg(S-\gamma_k(t))-arg(S-\gamma(t))\\
&\,\,\,\,-(arg(\gamma'(t))-arg(S-\gamma(t)))\\
=& ang(\gamma_k'(t),S-\gamma_k(t))-ang(\gamma'(t),S-\gamma(t))\\
&\,\,\,\,+arg(S-\gamma_k(t))-arg(S-\gamma(t))\,.
\end{align*}
By the previous part of the proof, we know that restrictions of $\gamma$ and $\gamma_k$ on $[a^k_j,b^k_j]$ are logarithmic in $V_{p_j}$ \textit{ if } $j\geq 2$.  Therefore, it follows (from the previous part of the proof) that
\begin{equation}\label{defdef}
Pr[ang(\gamma_k'(t),S-\gamma_k(t))-ang(\gamma'(t),S-\gamma(t))]=Pr[-\delta_{p_j}^k]\in (0,\delta),
\end{equation} 
if $j\geq 2\,$.
However, the same conclusion turns out to follow also if $j=1$, and thus
\begin{equation*}
Pr[ang(\gamma_k'(t),S_{V_0}-\gamma_k(t))-ang(\gamma'(t),S_{V_0}-\gamma(t))]\in (0,\delta),
\end{equation*}
if $t\in (a^k_1,b^k_1)=(w_k,b^k_1)$. 
This claim can be proven exatly by the same reasoning as we proved in the previous part of the proof that $-Pr[\delta_0^k]\in(0,\delta)$ (see also Picture 5). 

Then we are ready to prove (\ref{Crucial}). For this, let $\epsilon_k$ below denote generic $k$-dependent positive real numbers such that $\epsilon_k\to 0$ as $k\to\infty$. By the basic properties of the argument function, we get that 
\begin{equation*}
|Pr[arg(S-\gamma_k(t))-arg(S-\gamma(t))]|\leq \frac{1+\epsilon_k}{\min\{|S-\gamma_k(t)|,|S-\gamma(t)|\}}|\gamma_k(t)-\gamma(t)|\,. 
\end{equation*}
Because 
\begin{equation*}
\sup_{t\in[0,r]}\frac{|S-\gamma(t)|}{|S-\gamma_k(t)|}\to 1\,\text{ as }k\to\infty\,,
\end{equation*}
we conclude from above that
\begin{equation*}
Pr(ang(\gamma_k'(t),\gamma'(t)))\in (-\frac{(1+\epsilon_k)}{d(\gamma(t),\partial\Omega)}|\gamma_k(t)-\gamma(t)|\,,\,\delta\,)\,.
\end{equation*}

Recall then the formula used in the proof of Theorem \ref{ncase}, thus  
\begin{align*}
\frac{d}{dt}|\gamma_k(t)-\gamma(t)|&=\bigg(\frac{\gamma'_k(t)}{|\gamma_k'(t)|}-\frac{\gamma'(t)}{|\gamma'(t)|}\bigg)\cdot|\gamma_k'(t)|\frac{(\gamma_k(t)-\gamma(t))}{|\gamma_k(t)-\gamma(t)|}\\
&\,\,\,\,\,\,\,\,+\bigg(\frac{1}{|\gamma'(t)|}-\frac{1}{|\gamma_k'(t)|}\bigg)|\gamma_k'(t)|\,\frac{\gamma'(t)\cdot(\gamma_k(t)-\gamma(t))}{|\gamma_k(t)-\gamma(t)|}\,.
\end{align*}
Since $\big{|}|\gamma'(t)|-|\gamma_k'(t)|\big{|}\leq |\gamma_k(t)-\gamma(t)|$ and 
\begin{equation*}
 \frac{|\gamma'(t)\cdot(\gamma_k(t)-\gamma(t))|}{|\gamma_k(t)-\gamma(t)|}\leq \epsilon_k|\gamma'(t)|\,,
\end{equation*}
we conclude that the absolute value of the latter term can be estimated by
 $\epsilon_k|\gamma_k(t)-\gamma(t)|\,$. Furthermore, since 
\begin{equation*}
Pr[ang(\gamma'(t)\,,\,\gamma_k(t)-\gamma(t))]\approx -\frac{\pi}{2}\,,
\end{equation*}
it follows that 
\begin{align*}
\bigg(\frac{\gamma'_k(t)}{|\gamma_k'(t)|}-\frac{\gamma'(t)}{|\gamma'(t)|}\bigg)\cdot\frac{(\gamma_k(t)-\gamma(t))}{|\gamma_k(t)-\gamma(t)|}&\approx Pr[ang(\gamma_k'(t),\gamma'(t))]\\
&\geq -\frac{(1+\epsilon_k)}{d(\gamma(t),\partial\Omega)}|\gamma_k(t)-\gamma(t)|\,.
\end{align*}
Combining these estimates yields the desired inequality
if $t\in (a^k_j,b^k_j)$, $1\leq j\leq l$. This proves (\ref{Crucial}), since one can easily see that convergence $\gamma_k\to\gamma$ (and $\gamma'_k\to\gamma'$) implies that
\begin{equation}\label{loplop}
|\{\,t\in [w_k,r]\,:\,t\not\in\bigcup [a^k_j,b^k_j]\}|=:|A_k|\,\to 0\,,\text{ as }k\to\infty\,.
\end{equation}
\smallskip

Still we have to prove (\ref{Crucial2}), which corresponds to the case $t\in [b^k_j, a^k_{j+1}]$, $1\leq j\leq l-1$. Before this, we have to verify a couple of auxiliary results.

\subsubsection*{Auxiliary fact 1:} If $1\leq i\leq N$, then
\begin{equation}\label{claim1}
\gamma'(t_i) \text{ is not parallel to } \vec{J}_i,
\end{equation}
and, if $\gamma'(t_i)$ is parallel to $-\vec{J}_i$, then it holds for $k$ large enough that 
\begin{enumerate}
 \item $t_{i-1}<t_i$ and $\gamma$ is right-curving on $[t_{i-1},t_i]\,$,
\item  $t_i^k<t_i\,$ and $\gamma_k$ is left-curving on $(t_i^k,t^k_{i+1})\,$.
\end{enumerate}
For the proof of (\ref{claim1}), one suffices to verify that the counterclaim implies, by $Pr(\delta^k_i)>0$, that $\gamma'_k(t_i^k)$ actually points away from $V_{i+1}$. This clearly contradicts the above choice of $t_i,t_i^k,E_i$ and $V_i$ (according Lemma \ref{commoncells1}). The proofs of (1) and (2) follow easily from the definition of $t_i$ and assumptions (\ref{leftassumption}), (\ref{helpottava}) and (\ref{toka}). To become convinced of that, take a look at the illustrating pictures below.

\begin{figure}[htp]
\centering
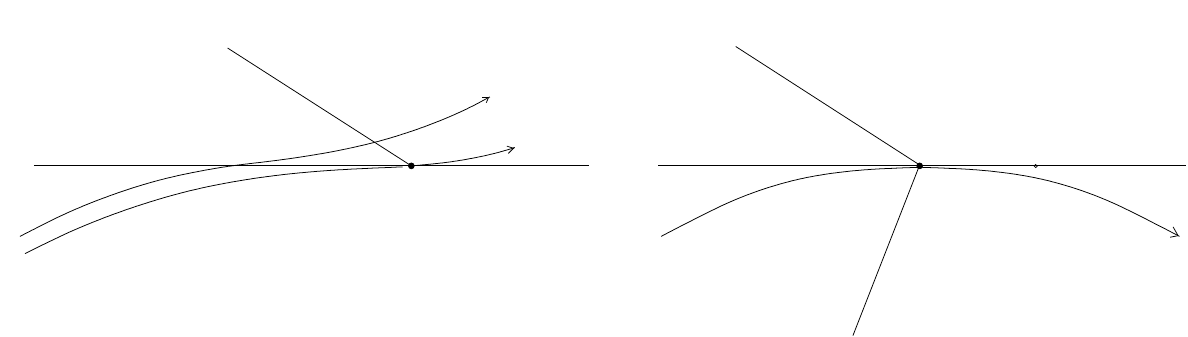
\label{kuva:10}
\end{figure}

Still, before the proof of (\ref{Crucial2}), we need the following result:
\subsubsection*{Auxiliary fact 2:} 
There exists a constant $C=C_{\Omega,\gamma}$ such that either
 \begin{equation}\label{option1}
a^k_{j+1}-b^k_j \leq C|\gamma(b^k_j)-\gamma_k(b^k_j)|\,
\end{equation}
or then $\gamma_k$ is left-curving and $\gamma$ right-curving on \[[b^k_j,\min\{t^{k}_{p_j+2},t_{p_j+1}\}]\,=:[b^k_j,\Lambda]\] and  
\begin{equation}\label{option2}
|a^k_{j+1}-\Lambda\}|\leq C|\gamma(\Lambda)-\gamma_k(\Lambda)|\,.
\end{equation}
This claim states that either the time the geodesics spend on different Voronoi cells (between the moments $b_j^k$ and $a^k_{j+1}$) is uniformly comparable to the distance of the geodesics at the moment of separation, or, then the geodesics curve away from each other some time (up to the moment $t=\Lambda$) after the separation, and then, before the geodesics again meet a common cell, the rest of the time of separation is again uniformly comparable to the distance at the moment of separation. The proof of this claim is given below, unless we remark that again the claim is more or less obvious in the light of the illustrating pictures below.
\begin{figure}[htp]
\centering
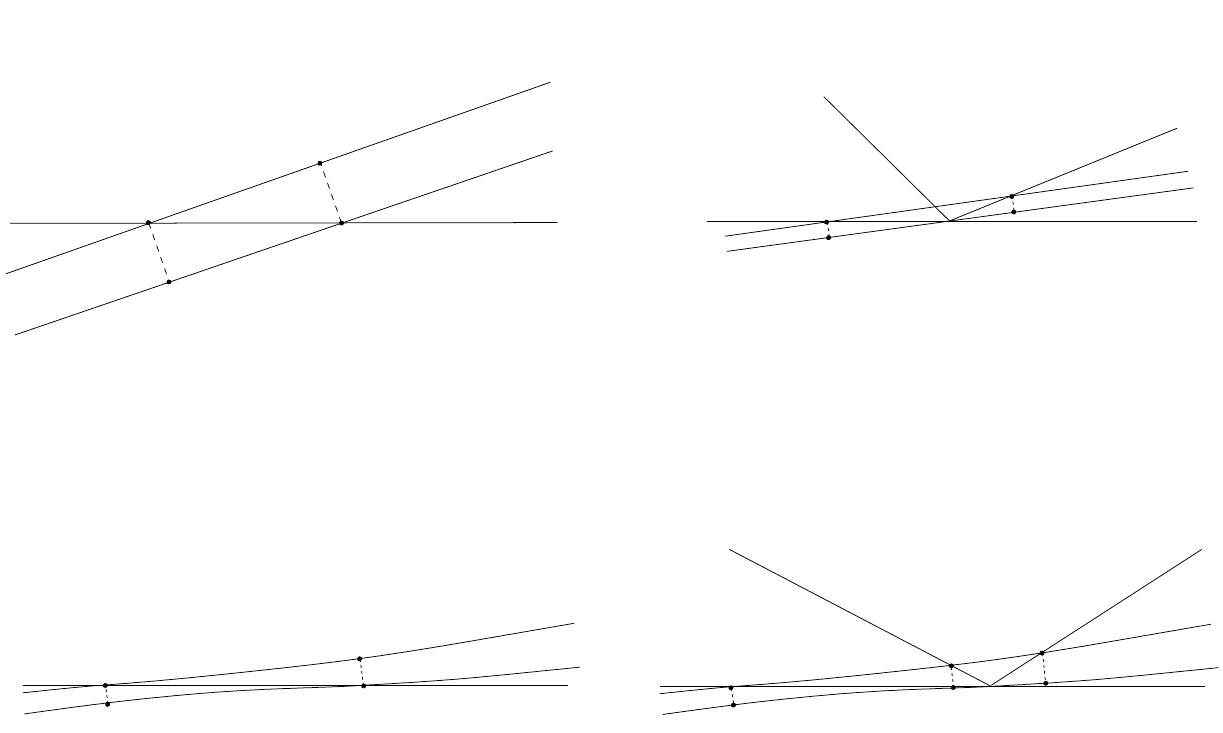
\label{kuva:10}
\end{figure} 
 
The first option (\ref{option1}) turns out to follow (as expected) if $\gamma'(t_{p_j+1})$ is not parallel to $-\vec{J}_{p_j+1}$. Therefore, option (\ref{option2}) corresponds to the case where $\gamma'(t_{p_j+1})$ is parallel to $-\vec{J}_{p_j+1}$, and thus we can use the results from above (auxiliary fact 1) in this case. 

To show the above claim, let below $C$ denote a generic constant which might depend on $\Omega$ and $\gamma$ but not on $k$. Consider first the claim for those $1\leq j\leq l-1$ for which $\gamma'(t_{p_j+1})$ is not parallel to $-\vec{J}_{p_j+1}=:-(Y_{p_j+1}-\gamma(t_{p_j+1}))$. By suitable translation and rotation, one can assume that $\gamma(t_{p_j+1})=0$, $\gamma'(t_{p_j+1})$ lies on the positive real axis, $|\gamma'(t_{p_j+1})|=d(\gamma(t_{p_j+1}),\partial\Omega)$, and $Im(Y_{p_j+1})>c_{\Omega,\gamma}>0$. 

Consider then the intersection point $P$ for a line $L(t)=\textbf{i}\delta_k+t$, $t\in \re$ (where $\delta_k>0$ converge to $ 0$ as $k\to 0\,$) and the line segment $J=[0,Y_{p_j+1}]$. Here $L$ approximates locally $\gamma_k$ and $P$ approximates the intersection of $\gamma_k$ and $J_{p_j+1}$. Clearly it holds that
$Re(P)\geq -C\delta_k\,$. By using again the assumptions (\ref{leftassumption}), (\ref{helpottava}) and (\ref{toka}) (as well as the fact that $\gamma$ and $\gamma_k$ are uniformly $C^1$) the above approximative sketch can be easily refined to show that 
\begin{equation}\label{latter2}
|t_{p_{j}+1}-b^k_j|\leq C\frac{|\gamma(b^k_j)-\gamma_k(b^k_j)|}{d(\gamma(t_{p_j+1}),\partial\Omega))}\,
\end{equation}
for $k$ large enough. Using (\ref{claim1}), exactly the same reasoning implies the same bound for $|a^k_{j+1}-t_{p_{j}+1}|$. Combining these implies (\ref{option1}).

For the case where $\gamma'(t_{p_j+1})$ is parallel to $\vec{J}_{p_j+1}=-(Y_{p_j+1}-\gamma(t_{p_j+1}))$ it follows directly from the auxiliary fact 1 that $b_j^k=\min\{t_{p_j+1}, t^k_{p_j+1}\}=t^k_{p_j+1}\,$, $\gamma_k$ is left-curving and $\gamma$ is right-curving on $[b_j^k,\Lambda]$. Furthermore, (\ref{option2}) follows essentially by the same argument as was used to prove (\ref{option1}) in the case of non-parallelity. More precisely, notice that if $\Lambda=t_{p_j+1}$, then the claim equals to the case of $|a^k_{j+1}-t_{p_{j}+1}|$ above. In turn, if $\Lambda=t^k_{p_j+2}$ (for large $k$), then it is easy to see that $\gamma(t_{p_j+1})$ must lie in a corner, $t_{p_j+2}=t_{p_j+1}$ and $\gamma'(t_{p_j+2})$ is not parallel to $-\vec{J}_{p_j+2}$.
Then one only has to repeat the argument for (\ref{option1}) from above to complete the proof.  
\smallskip

Before the final part of the proof, let us denote
\begin{equation*}
G(t)=|\gamma_k(t)-\gamma(t)|\,.
\end{equation*}
By the uniform Lipschitz estimate (Proposition \ref{apu2}) for $\gamma'$ and $\gamma_k'$ it is easy to check that
\begin{equation}\label{apu99}
|G'(t_1)-G'(t_2)|\leq C|t_1-t_2|\,,\text{ if }t_1,t_2\in[0,r]\,.
\end{equation}


\subsubsection*{The proof of (\ref{Crucial2})}
 
Suppose first that $t\in [b^k_j,a^k_{j+1}]$ such that (\ref{option1}) holds. In this case we have (recall the generic choice of $C$)
\begin{align*}
G'(t)&\geq G'(b_j^k)-|G'(t)-G'(b_j^k)|\geq G'(b_j^k)-C|t-b_ j^k|\\
&\geq G'(b_j^k)-CCG(b_j^k)\geq -(1+\epsilon_k)G(b_j^k)-CG(b_j^k)
\geq -CG(b_j^k)\,.
\end{align*}
This implies that
\begin{align*}
G(t)=G(b_j^k)+\int_{b_j^k}^tG'(t)\,dt\,\geq G(b_j^k)-|t-b_j^k|CG(b_j^k)\geq G(b_j^k)(1-\epsilon_k)\,.
\end{align*}
The last inequality above follows by $|a^k_{j+1}-b^k_j|\to 0$ as $k\to\infty$. Finally, applying this inequality again back to the first inequality we get that 
\begin{equation*}
G'(t)\geq C(-G(b_j^k))\geq C(-\frac{1}{1-\epsilon_k}G(t))\geq -CG(t)\,, 
\end{equation*}
which proves the claim. Yet we have to show the claim in the case of option (\ref{option2}). For this, suppose first that $t\leq \Lambda$. Recall from above that in the case of option (\ref{option2}), $\gamma'(t_{p_j+1})$ is parallel to $-\vec{J}_{p_j+1}$, and $t\leq\Lambda$ implies $t\in[t_{p_j}, t_{p_j+1}]$. Therefore, if $S$ denotes the nucleus of $V_{p_j}$, we may write
\begin{align*}
&ang(\gamma_k'(t),S-\gamma_k(t))-ang(\gamma'(t),S-\gamma(t))\\
=&ang(\gamma_k'(t),S-\gamma_k(t))-ang(\gamma'(b_j^k),S-\gamma(b_j^k))\\
=&ang(\gamma_k'(b_j^k),S-\gamma_k(b_j^k))-ang(\gamma'(b_j^k),S-\gamma(b_j^k))\\
&\,\,\,\,+ang(\gamma_k'(t),S-\gamma_k(t))-ang(\gamma_k'(b_j^k),S-\gamma_k(b_j^k))\\
=&-\delta^k_{p_j}+ang(\gamma_k'(t),S-\gamma_k(t))-ang(\gamma_k'(b_j^k),S-\gamma_k(b_j^k))\,\\
=&-\delta^k_{p_j}+ang(\gamma_k'(t),\gamma_k'(b_j^k))+ang(S-\gamma_k(b_j^k),S-\gamma_k(t))\,.
\end{align*}
We proceed with showing that the principal value of the above quantity lies in $(0,\delta)$. For this, we have shown above that $Pr(-\delta^k_{p_j})\in (0,\delta)\,$. Moreover, since $\gamma_k$ is left-curving on $[b^k_j,\Lambda]$ (by auxiliary fact 2) and $t>b^k_j$, it holds that $Pr[ang(\gamma_k'(t),\gamma_k'(b_j^k))]\in (0,\delta)$. 

To show that 
\begin{equation}\label{tikka}
Pr[ang(S-\gamma_k(b_j^k),S-\gamma_k(t))\,]\in(0,\delta)\,, 
\end{equation}
observe that by a suitable rotation and translation, we may assume that $J_{p_j+1}$ is contained in the positive real-axis, $\gamma(t_{p_j+1})=0$, $V_{p_j}$ is contained in the upper half-plane $\{z:Im(z)\geq 0\}$ and, especially, $Im(S)>c>0\,$, where $c$ is independent of $k$.
Since $\gamma_k(b_j^k)\in J_{p_j+1}$ lies on the positive real-axis and converge to $\gamma(t_{p_j+1})=0$ as $k\to\infty$, it follows that $Im(S-\gamma_k(b_j^k))>c$ for $k$ large enough. Combining this with the fact that 
\begin{equation*}
Pr[\gamma_k(b^k_j)-\gamma_k(t)]\to 0\,\text{ as }k\to\infty\,,  
\end{equation*}
clearly implies that 
\begin{align*}
&Pr[ang(S-\gamma_k(b_j^k),S-\gamma_k(t))]\\
=&Pr[arg\big(S-\gamma_k(b_j^k)\big)-arg\big(S-\gamma_k(b_j^k)+(\gamma_k(b_j^k)-\gamma_k(t))\big)]\in(0,\delta)\,.
\end{align*}
Thus, we have shown that 
\begin{equation*}
Pr[ang(\gamma_k'(t),S-\gamma_k(t))-ang(\gamma'(t),S-\gamma(t))]\in(0,\delta)\,.
\end{equation*}
Then exactly the same reasoning as in the beginning of the proof of (\ref{Crucial}) gives that
\begin{equation*}
Pr(ang(\gamma_k'(t),\gamma'(t)))\in\big(\,-\frac{1+\epsilon_k}{d(\gamma(t),\partial\Omega)}|\gamma_k(t)-\gamma(t)|,\,\delta\,\big)\,
\end{equation*}
if $t\in[b_j^k,\Lambda]\,$ and also the inequality (\ref{Crucial}) (which implies (\ref{Crucial2})) if $t\in[b_j^k,\Lambda]\,$. 

The remaining case $t\in[\Lambda,a^k_{j+1}]$ then turns out to follow essentially in the same way as the case of (\ref{option1}) above. Indeed, first recall that by auxiliary fact 2, it holds that $|a^k_j-\Lambda|\leq C|\gamma_k(\Lambda)-\gamma(\Lambda)|\,$ and, as it was shown just above, 
\begin{equation*}
G'(\Lambda)\geq -(1+\epsilon_k)G(\Lambda)\,,
\end{equation*}
where $G(t)=|\gamma_k(t)-\gamma(t)|\,$, as we defined above. Then one can repeat the argument for the case of (\ref{option1}) above (replace $b_j^k$ by $\Lambda$ in the formulas) to obtain (\ref{Crucial2}). This completes the proof. 
\hfill$\Box$\\

When utilizing the previous lemma, we are going to exploit the following version of the classical Gronwall's inequality(\cite{Gr}):
\begin{proposition}\label{gronwall}
Suppose that $G:[0,r]\to [0,\infty)$ is $C^1$, $G(0)=0$, and $g:[0,r]\to\re$ is integrable such that
\begin{equation*}
G'(t)\geq g(t)G(t)\,\text{ for all }t\in[0,r].  
\end{equation*}
Then 
\begin{equation*}
G(r)\geq G(t)\,exp\bigg(\int_{t}^{r}g(z)\,dz\,\bigg)\,\,,\text{ for all }t\in[0,r]\,.
\end{equation*} 
\end{proposition}

\begin{corollary}\label{major3}
 Under the conditions of Lemma \ref{major}, it follows that
\begin{equation*}
\max_{t\in[0,r]}|\gamma_k(t)-\gamma(t)|\leq e^{(1+\epsilon_k)r}|y-y_k|\,,
\end{equation*}
where $\epsilon_k\to 0$ as $k\to\infty\,$.
\end{corollary}
\textit{Proof.}
The claim follows from Lemma \ref{major} by applying Proposition \ref{gronwall} with
\begin{equation*}
G(t)=|\gamma_k(t)-\gamma(t)|\,\,\text{ and }\,g(t)=-(1+\epsilon_k)-\chi_{A_k}(t)C_{\gamma,\Omega}\,.
\end{equation*}
\hfill$\Box$\\

\begin{lemma}\label{tili}
Suppose that $\Omega$ is a domain in $\re^2$ (or in $\rn$) and let the canonically parametrized geodesics $\gamma:[0,r]\to\Omega$ from $x$ to $y$, 
 and $\tilde{\gamma}:[0,r]\to\Omega$ from $x$ to $\tilde{y}\neq y$ satisfy
\begin{equation}\label{ekaekaeka}
\max_{t\in[0,r]}\big|\gamma(t)-\tilde{\gamma}(t)\big|\leq C|y-\tilde{y}|\,,
\end{equation}
and 
\begin{equation}\label{tokatoka}
|y-\tilde{y}|\leq \frac{1}{C}\inf_{z\in B_Q(x,r)}d\big(z,\partial\Omega\big)\,=:\frac{1}{C}m\,. 
\end{equation}
Then 
\begin{equation*}
l_Q(\frac{\gamma(t)+\tilde{\gamma}(t)}{2})\leq r+\frac{rC^2}{m^2}|y-\tilde{y}|^2\,.
\end{equation*}
\end{lemma}
\textit{Proof.} Let $\omega$ denote the average path of $\gamma$ and $\tilde{\gamma}$.  
It is easy to check that by (\ref{ekaekaeka}) and
(\ref{tokatoka}) we have
\begin{equation*}
d(\omega(t),\partial\Omega)\geq \frac{m}{2}\,\,\text{ for all }t\in[0,r]\,.
\end{equation*}
Notice also that 
\begin{equation*}
\,\,d(\gamma(t),\partial\Omega)\geq m\,\text{ and }\,d(\tilde{\gamma}(t),\partial\Omega)\geq m\,\text{ if }t\in[0,r]
\end{equation*}
by the assumption for $m\,$.
By Proposition \ref{distcurv}, we have for all $t\in[0,r]$ that
\begin{align*}
d\big(\frac{\gamma(t)+\tilde{\gamma}(t)}{2},\partial\Omega\big)&\geq
\frac{d(\gamma(t),\partial\Omega)+d(\tilde{\gamma}(t),\partial\Omega)}{2}-\frac{\big(\frac{|\gamma(t)-\tilde{\gamma}(t)|}{2}\big)^2}{2d(\frac{\gamma(t)+\tilde{\gamma}(t)}{2},\partial\Omega)}\\
&\geq\frac{|\gamma'(t)|+|\tilde{\gamma}'(t)|}{2}-\frac{|\gamma(t)-\tilde{\gamma}(t)|^2}{4m}\,. 
\end{align*}
This yields the claim by
\begin{align*}
 l_Q(\omega)&\leq \int_{0}^{r}\frac{|\gamma'(t)|+|\tilde{\gamma}'(t)|}
{2d\big(\frac{\gamma(t)+\tilde{\gamma}(t)}{2},\partial\Omega\big)}\,dt\\ 
&\leq 
\int_{0}^{r} 1+\frac{|\gamma(t)-\tilde{\gamma}(t)|^2}
{2m\big(|\gamma'(t)|+|\tilde{\gamma}'(t)|\big)-|\gamma(t)-\tilde{\gamma}(t)|^2}\,dt\\
&\leq r+\int_{0}^{r}\frac{|\gamma(t)-\tilde{\gamma}(t)|^2}{3m^2}\,dt\,\\
&\leq r+\frac{rC^2}{m^2}|y-\tilde{y}|^2\,.
\end{align*}
\hfill$\Box$\\

\begin{lemma}\label{versio}
Suppose that $\Omega$ is a Voronoi domain and there exists $x,y\in\Omega$ with two disjoint geodesics $\gamma_1$ and $\tilde{\gamma}_1$ from $x$ to $y$. Then $D_{\gamma_1\cup\tilde{\gamma}_1}$ contains $z\in\partial\Omega\,$. 
\end{lemma}
\textit{Proof.} Suppose that the claim is not true. Then the conditions are exactly as in Lemma \ref{topologia2}, and thus we can directly apply it and conclude with $\gamma$ and sequences $\gamma_k$ and $\tilde{\gamma}_k$ according to the statement of the lemma.
Let $r$ denote the quasihyperbolic length of $\gamma$ and $r_k$ the quasihyperbolic length of $\gamma_k$ and $\tilde{\gamma}_k$. Since 
\begin{equation}\label{loop}
\gamma\subset \bar{D}_{\gamma_k\cup\tilde{\gamma}_k}\text{ for all }k\in\na\,,
\end{equation}
and $\gamma_k,\tilde{\gamma}_k$ are uniformly $C^1$, it is easy to check that $r_k\geq r$ for large $k$. Actually, the case $r_k=r$ can occur (for large $k$) only if either $\gamma_k$ or $\tilde{\gamma}_k$ coincide with $\gamma\,$. As in Lemma \ref{topologia2}, denote by $t_k$ the points of divergence for $\gamma_k$ and $\tilde{\gamma}_k\,$, and choose $t_k<\tilde{r}<r\,$ arbitrarily. 

By Lemma \ref{major3}, we can estimate that
\begin{align*}
\max_{t\in[0,\tilde{r}]}|\gamma_k(t)-\tilde{\gamma}_k(t)|&\leq \max_{t\in[0,\tilde{r}]}|\gamma_k(t)-\gamma(t)|+\max_{t\in[0,\tilde{r}]}|\tilde{\gamma}_k(t)-\gamma(t)|\\
&\leq e^{2\tilde{r}}\big(|\gamma_k(\tilde{r})-\gamma(\tilde{r})|+|\tilde{\gamma}_k(\tilde{r})-\gamma(\tilde{r})|\big)\,
\end{align*}
for all $k$ large enough. Then, if $\gamma_k(\tilde{r})=\gamma(\tilde{r})$ or $\tilde{\gamma}_k(\tilde{r})=\gamma(\tilde{r})$, it directly follows that
\begin{equation}\label{uiui}
\max_{t\in[0,\tilde{r}]}|\gamma_k(t)-\tilde{\gamma}_k(t)|\leq
e^{2\tilde{r}}\big|\gamma_k(\tilde{r})-\tilde{\gamma}_k(\tilde{r})\big|\,.
\end{equation}
On the other hand, if  $\gamma_k(\tilde{r})\neq\gamma(\tilde{r})\neq\tilde{\gamma}_k(\tilde{r})$, essentially the same conclusion turns out to follow. To see this, notice that by (\ref{loop}) $\gamma_k(\tilde{r})$ and $\tilde{\gamma_k}(\tilde{r})$ lie on the 'different sides of $\gamma$'. More precisely, we may assume, by suitable translation and rotation, that $\gamma(\tilde{r})=0$ and $\gamma'(\tilde{r})=c\textbf{i}$, $c>0$. Then, by (\ref{loop}), we have 
\begin{equation}\label{horohoro}
\frac{1}{2\pi}\bigg|\int_{\gamma_k\cup\tilde{\gamma}_k}\frac{dz}{z}\,\bigg|\,=1\,.
\end{equation}
By combining this with the fact that $\gamma_k,\tilde{\gamma}_k\to \gamma$ and $\gamma'_k,\tilde{\gamma}'_k\to \gamma'$ uniformly as $k\to\infty\,$,
it is easy to see that for $k$ large enough (\ref{horohoro}) is possible only if 
\begin{equation}
\frac{\gamma_k(\tilde{r})}{|\gamma_k(\tilde{r})|}\approx \pm 1\,\text{ and }\,
\frac{\tilde{\gamma}_k(\tilde{r})}{|\tilde{\gamma}_k(\tilde{r})|}\approx \mp 1\,.
\end{equation}
This implies the desired conclusion, thus
\begin{align*}
e^{2\tilde{r}}\big(|\gamma_k(\tilde{r})-\gamma(\tilde{r})|+|\tilde{\gamma}_k(\tilde{r})-\gamma(\tilde{r})|\big)\,
\leq e^{2\tilde{r}}2\big|\gamma_k(\tilde{r})-\tilde{\gamma}_k(\tilde{r})\big|\, 
\end{align*}
for all $k$ large enough. Summing up, we have shown that for all $k$ large enough the assumption (\ref{ekaekaeka}) in Lemma \ref{tili} holds with $C=2e^{2r}\,$.
Furthermore, let 
\begin{equation*}
 m:=\inf_{y\in B_Q(x,r)}d(y,\partial\Omega)\,\text{ and }\omega_k(t)=\frac{\gamma_k(t)+\tilde{\gamma}_k(t)}{2}\,,\,\,\,t\in[0,\tilde{r}]\,,
\end{equation*}
and choose $k_0$ to be so big that 
\begin{equation*}
|\gamma_k(\tilde{r})-\tilde{\gamma}_k(\tilde{r})|\leq \frac{m}{2e^{2r}}\,
\end{equation*}
for all $k\geq k_0\,$. Then we can directly apply Lemma \ref{tili} to obtain
\begin{equation*}
l_Q(\omega_k)\leq \tilde{r}+\frac{4re^{4r}}{m^2}|\gamma_k(\tilde{r})-\tilde{\gamma}_k(\tilde{r})|^2\, 
\end{equation*}
for all $k$ large enough.

For the desired contradiction, let us denote $M:=\sup_{w\in B_Q(y,r_k-\tilde{r})}d(w,\partial\Omega)$ and recall Theorem \ref{smallballs}, revealing that
\begin{equation*}
d_{Q}\bigg(\frac{\gamma_k(\tilde{r})+\tilde{\gamma}_k(\tilde{r})}{2},\gamma_k(r_k)\bigg)
\leq (r_k-\tilde{r})-\frac{|\gamma_k(\tilde{r})-\tilde{\gamma}_k(\tilde{r})|^2}{512(r_k-\tilde{r})M^2}\,.
\end{equation*}
Combining the above estimates, we get that 
\begin{align*}
r_k&=d_{Q}(x,\gamma_k(r_k))\leq d_Q(x,\omega_k(\tilde{r}))+d_Q(\omega_k(\tilde{r}),\gamma_k(r_k))\\ &\leq\,\tilde{r}+\frac{4re^{4r}}{m^2}|\gamma_k(\tilde{r})-\tilde{\gamma}_k(\tilde{r})|^2+(r_k-\tilde{r})-\frac{|\gamma_k(\tilde{r})-\tilde{\gamma}_k(\tilde{r})|^2}{512(r_k-\tilde{r})M^2}\\
&\leq \,r_k+|\gamma_k(\tilde{r})-\tilde{\gamma}_k(\tilde{r})|^2\bigg(\frac{4re^{4r}}{m^2}-\frac{1}{512(r_k-\tilde{r})M^2}\bigg)\,.
\end{align*}
Recalling that $r_k\to r$ as $k\to\infty$ and $r-\tilde{r}$ can be chosen as small as desired, it follows that
\begin{equation*}
\frac{4re^{4r}}{m^2}-\frac{1}{512(r_k-\tilde{r})M^2}<0\,,
\end{equation*}
implying that $r_k=d_{Q}(x,\gamma_k(r_k))<r_k$. This is the desired contradiction. 
\hfill$\Box$

\begin{corollary}\label{Voronoissa1}
If $\Omega$ is a Voronoi domain and $d_Q(x,y)<\pi$, then there exists exactly one geodesic from $x$ to $y$.  
\end{corollary}
\textit{Proof.}
On the contrary, suppose that there exists two different geodesics $\gamma:[0,r]\to\Omega$ and $\tilde{\gamma}:[0,r]\to\Omega$ from $x$ to $y$ such that $r<\pi$ ($\gamma,\tilde{\gamma}$ are canonocally parametrized). Since $\gamma\neq\tilde{\gamma}$, there exists $0\leq r_1<r_2\leq r$ such that $\gamma(t)\neq\tilde{\gamma}(t)$ if $t\in(r_1,r_2)$, $\gamma(r_1)=\tilde{\gamma}(r_1)$ and $\gamma(r_2)=\tilde{\gamma}(r_2)$. Let $\gamma_{|[r_1,r_2]}=:\omega$ and $\tilde{\gamma}_{|[r_1,r_2]}=:\tilde{\omega}$. Then the previous lemma implies that 
\begin{equation}\label{circle}
D_{\omega\cup\tilde{\omega}}\cap\partial\Omega\,\neq\,\emptyset\,.
\end{equation}
This yields the desired contradiction, since it is well known that (\ref{circle}) implies that 
\begin{equation*}
l_Q(\omega\cup\tilde{\omega})\geq 2\pi>2r\geq 2(r_2-r_1)=l_Q(\omega\cup\tilde{\omega})\,.
\end{equation*}
\hfill$\Box$\\

For the following corollary concerning the case of simply connected domains, we make a convention that notation $\Omega_k\to\Omega\subset\re^2$ means (in this work) that $\Omega_k$ is a decreasing sequence of Voronoi domains such that $\partial\Omega_k\subset\partial\Omega$ for every $k\in\na$ and for every $R>0$ it holds that
\begin{equation*}
\sup_{y\in \partial\Omega\cap B(0,R)}d(y,\partial\Omega_k)\,\to 0\text{ as }k\to\infty\,.
\end{equation*}

\begin{lemma}\label{Voronoissa}
 Suppose that $\Omega\subset\re^2$ is a simply connected domain, $\Omega_k\to \Omega$,  $x\in\Omega$, $r>0$ and let $B_k(x,r)$ denote the quasihyperbolic ball related to $\Omega_k$. Then there exists $k_0\in\na$ such that if $k\geq k_0$, then for every $y\in B_k(x,r)$ there exists a unique geodesic from $x$ to $y$ (with respect to the quasihyperbolic metric in $\Omega_k$). 
\end{lemma}
\textit{Proof.} Let $x\in\Omega$, $r>0$, $k\in\na$. First notice that there exists $C,c>0$ such that 
\begin{equation}\label{apu1}
M_k:=\sup_{y\in B_k(x,r)}|y|\,<\, C\,\text{ and }\,m_k:=\inf_{y\in B_k(x,r)}d(y,\partial\Omega_k)>c>0\,.
\end{equation}
Both of these claims are well known. For example, to see the claim for $m_k$, observe that $x\in\Omega$ implies $d(x,\partial\Omega_k)>c'>0$ for all $k$, and, if (\ref{apu1}) does not hold, then there exist $a_k\in B_k(x,r)$ such that $d(a_k, \partial\Omega_k)\to 0$ as $k\to\infty$. It is well known that in this case $d_{Q}(x,a_k)$ (with respect to $\Omega_k$) tends to infinity, which contradicts with $a_k\in B_k(x,r)$. 

Suppose then that $B_k(x,r)\subset B(0,C)$. Let us then choose $k_0$ to be so large that the bounds in (\ref{apu1}) are valid for all $k\geq k_0$ and for all $z\in \partial\Omega\cap B(0,C)$ there exists $z_k\in \partial\Omega_k$ such that 
\begin{equation}\label{distdist}
|z-z_k|\leq\frac{c}{2}\,.
\end{equation}

Let then $k\geq k_0$ and fix a point $y\in B_k(x,r)\,$. Suppose that there exist two different geodesics $\gamma_k$ and $\tilde{\gamma}_k$ between $x$ and $y$. As in the previous corollary, this yields the existence of $\omega_k$ and $\tilde{\omega}_k$ such that $\omega_k$ and $\tilde{\omega}_k$ are disjoint, $\omega_k\cup\tilde{\omega}_k\subset B_k(x,r)$, and there exists $a\in D_{\omega_k\cup\tilde{\omega}_k}\cap\partial\Omega_k\,$. Since $\partial\Omega_k\subset \partial\Omega$, we thus get that $a\in D_{\omega_k\cup\tilde{\omega}_k}\cap\partial\Omega\,$. Since $\Omega$ is simply connected (and $x\in\Omega\cap(\omega\cup\tilde{\omega})$), we conclude that $\omega_k\cup\tilde{\omega}_k$ has to intersect the boundary of $\Omega$, and thus there exists $z\in B_k(x,r)\cap\partial\Omega$. By (\ref{distdist}), there exists $z_k\in\partial\Omega_k$ such that $|z-z_k|\leq \frac{c}{2}\,$. This contradicts with (\ref{apu1}) and completes the proof.
\hfill$\Box$\\


\section{The case of general domains}\label{generaldomains}
For the proof of conjectures \ref{konj2} and \ref{konj3}, the main work has been done in the previous section but still some auxiliary lemmata are needed. To be able to extend the results from the previous section to the case of general domains, we have to extend Corollary \ref{major3} from the local scale to the 'semi-global' scale. More precisely, we will consider the situation where $\Omega$ is a domain in $\re^2$, $x\in\Omega$ and $r+1>\tilde{r}>r>0$ so that for every $y\in B_Q(x,\tilde{r})$ there exists a unique geodesic $\gamma_{x,y}$. Then we seek for an appropriate quantity $c>0$, so that if $y,z\in\partial B_Q(x,r)$ such that $|y-z|\leq c$, then there exists a $2$-Lipschitz curve from $[0,|z-y|]$ to $\partial B_Q(x,r)$ connecting $z$ and $y$. 

For the forthcoming results, let us introduce the following notation: if there exists a unique quasihyperbolic geodesic $\gamma$ between $x$ and $y$ in $\Omega$, then we denote 
\begin{equation*}
\gamma=:\gamma_{x,y}\,.
\end{equation*}

Also denote, if $\Omega$ is a domain in $\re^2$, $x\in\Omega$, $0<r<\tilde{r}<r+1$, that
\begin{equation*}
m(\Omega,x,\tilde{r})=\inf_{w\in B_{Q}(x,\tilde{r})}d(w,\partial\Omega)\,,
\end{equation*}
and
\begin{equation*}
c(\Omega,x,r,\tilde{r}):=\frac{\min\big{\{}\,1\,,\,\big(m(\Omega,x,\tilde{r})\big)^3\big{\}}(\tilde{r}-r)}{10^{10}\max\{\,1\,,\,4e^{2r}\,\}}\,.
\end{equation*}
The main idea in this final part is to deduce the following lemma from Corollary \ref{major3}:
\begin{lemma}\label{loppusjo}
Suppose that $\Omega$ is a Voronoi domain in $\re^2$, $x\in\Omega$ and $r+1>\tilde{r}>r>0$ such that for every $y\in B_Q(x,\tilde{r})$ there exists $\gamma_{x,y}$. Then, if
\begin{equation*}
y,z\in \partial B_Q(x,r)\,\text{ such that }\,|y-z|\leq c(\Omega,x,r,\tilde{r}),
\end{equation*}
then   
\begin{equation*}
d_Q(x,\frac{y+z}{2})\leq \,r+\frac{4re^{4r}}{(m(\Omega,x,\tilde{r}))^2}|y-z|^2\,.
\end{equation*}
\end{lemma}
This result, which is proved later, guarantees an appropriate bound for the anticonvexity of quasihyperbolic balls with unique geodesics in a slightly bigger ball. The proof of the lemma turns out to follow rather easily from the results of the previous section and the following result:


\begin{lemma}\label{loppusjo3}
Suppose that $\Omega$ is a Voronoi domain, $x\in\Omega$ and $B_Q(x,r)$ is such that there exists $r+1>\tilde{r}>r>0$ such that for every $y\in B_Q(x,\tilde{r})$ there exists a unique geodesic $\gamma_{x,y}$. Then  
\begin{equation*}
\max_{t\in[0,r]}|\gamma_{x,y}(t)-\gamma_{x,z}(t)|\leq\,2e^{2r}|y-z|\,,
\end{equation*}
for all $y,z\in\partial B_Q(x,r)\,$ such that $|y-z|\leq c(\Omega,x,r,\tilde{r})\,$.
\end{lemma}

Lemma \ref{loppusjo3} is a semi-global version of Lemma \ref{major} and its corollary \ref{major3}. To prove it, we first have to verify some elementary estimates related to the chosen quantity $c(\Omega,x,r,\tilde{r})$. We remark that the used arguments are obviously very rough. They might also be more or less well known but we did not manage to find suitable references to be able to avoid proving them. However, the details of the proofs are partially left to the reader. We also remark that the corresponding (and more sharp) analysis can be carried out in general dimensions.




\begin{proposition}
Suppose that $\Omega$ is a domain in $\re^2$ and $x,y\in\Omega$, $d_Q(x,y)=r>0$, are such that there exists $\varepsilon>0$ such that the unique geodesic $\gamma_{x,\tilde{y}}$ exists for all $\tilde{y}\in B_Q(y,\varepsilon)$. Then there exists a prolongation of $\gamma_{x,y}:[0,r]\to\Omega$, that is a geodesic $\gamma:[0,\tilde{r}]\to\Omega$, $l_Q(\gamma)=\tilde{r}>r$,  such that $\gamma(t)=\gamma_{x,y}(t)$ if $t\in[0,r]\,$.
\end{proposition}
\textit{Sketch of the proof.}
Let us denote that $\gamma_{x,y}=:\gamma\,$. First of all, if $\gamma_k$ is a sequence of geodesics from $x$ to $y_k$ and $y_k\to y$ as $k\to\infty$, then the assumptions imply that $\gamma_k\to\gamma$ and also that $\gamma'_k\to\gamma'$ as $k\to\infty$. Therefore, if $\epsilon>0$, we can choose $\delta>0$ such that if $\tilde{y}\in B_Q(y,\delta)\,$, it follows that
\begin{equation}
\sup_{t\in [r-\delta,r+\delta]}\big|\gamma_{\tilde{y}}'(t)-\gamma'(r)\big|\,\leq\epsilon\,.
\end{equation}
Furthermore, 
we may assume that $\gamma(r)=y=0$ and $\frac{\gamma'(r)}{|\gamma'(r)|}=\textbf{i}$. Consider then the line 
\begin{equation*}
L(t)=\delta_1\gamma'(r)+t\,,\,\,\,t\in\re\,,\,\delta_1<<\varepsilon\,.
\end{equation*}
It is easy to check that for every $t\in [-\delta_1,\delta_1]$  there exists a geodesic $\gamma_t$ from $x$ to $L(t)$ of length strictly greater than $r$. Moreover, it is clear
that for every $t\in[-\delta_1,\delta_1]$, $\gamma_t$ intersects the real-axis, at a point denoted by $(a(t),0)$ and such that $a(-\delta_1)< 0$ and $a(\delta_1)> 0$. Suppose then that the claim is not true, implying that $|a(t)|>c>0$ for all $t\in[-\delta_1,\delta_1]\,$. 
Then by defining
\begin{equation*}
\tilde{t}=\inf\{t\in [-\delta_1,\delta_1]\,:\,a(t)>0\}\,,
\end{equation*}
one can easily verify the existence of two disjoint geodesics $\omega$ and $\tilde{\omega}$ from $x$ to $L(\tilde{t})$. This contradicts the assumption and completes the proof.  
\hfill$\Box$\\

\begin{corollary}
Suppose that $\Omega$ is a domain in $\re^2$, $x\in\Omega$ and $r>0$ so that for every $y\in B_Q(x,r)$ there exists a unique geodesic $\gamma_{x,y}$. Then, for every $y\in B_Q(x,r)$, $\gamma_{x,y}$ is a subarc of a quasihyperbolic geodesic $\gamma_{x,y'}$ of length $r$.
\end{corollary}
\textit{Proof.} The above corollary implies that each geodesic $\gamma_{x,y}$ with $l_Q(\gamma_{x,y})<r$ has a prolongation. To see that there exists a prolongation of length $r$, this result can be iterated to obtain an increasing sequence of prolongations $\gamma_k$ such that at every step the length of the prolongations tends closer and closer to the maximal length of the prolongations. Then it is easy to check that $(\gamma_k)$ converges and $l_Q(\gamma)=r$ for the limit $\gamma$.
\hfill$\Box$\\

\begin{proposition}\label{eiviitetta}
Suppose that $\Omega$ is a domain in $\re^2$, $x\in\Omega$ and $\tilde{r}>r>0$ so that for every $y\in B_Q(x,\tilde{r})$ there exists a unique geodesic $\gamma_{x,y}$. Then, if $y,z\in\partial B_Q(x,r)$  such that 
\begin{equation}\label{ddoletus}
|y-z|\leq 2c(\Omega,x,r,\tilde{r})\,,
\end{equation}
it holds that
\begin{equation}\label{aaclaim}
\bigg|\frac{\gamma_{x,y}'(r)}{|\gamma_{x,y}'(r)|}\cdot\frac{(y-z)}{|y-z|}\bigg|\leq\frac{1}{100}\,.
\end{equation}
\end{proposition}
\textit{Proof.}
By a suitable translation, we may assume that $y=0$. Suppose, on the contrary that $y,z\in\partial B_Q(x,r)$ satisfy (\ref{ddoletus}) but (\ref{aaclaim}) does not hold. Assume first that 
\begin{equation}\label{ekakeissi}
\frac{\gamma_{x,y}'(r)}{|\gamma_{x,y}'(r)|}\cdot\frac{z}{|z|}>\frac{1}{100}\,.
\end{equation}
Also denote, for simplicity, that $\gamma_{x,y}=:\gamma$, let $\tilde{\gamma}$ be the prolongation of $\gamma$ with quasihyperbolic length $\tilde{r}$ and 
\begin{equation*}
\tilde{t}:=\frac{2000|z|}{d(0,\partial\Omega)}\,\,\text{ and }\,a=\tilde{\gamma}(r+\tilde{t})\,.
\end{equation*}
Assumption (\ref{ddoletus}) guarantees (recall that $y=0$) that $r+\tilde{t}<\tilde{r}\,$, as well as $\tilde{t}<<1$.
Next we are going to show that (\ref{ddoletus}) combined with (\ref{ekakeissi}) implies that 
$d_Q(z,a)<\tilde{t}\,$, which yields the desired contradiction. To show this, observe first that (by Gehring-Palka)
\begin{equation*}
|a|\leq (e^{\tilde{t}}-1)d(0,\partial\Omega)<\tilde{t}e^{\tilde{t}}d(0,\partial\Omega)\,.
\end{equation*}
This implies, especially, that $J[0,a]\subset \Omega$, and also that  
\begin{align*}
&\frac{|a|}{(1+\tilde{t}e^{\tilde{t}})d(0,\partial\Omega)}\leq\tilde{t}\leq \frac{|a|}{(1-\tilde{t}e^{\tilde{t}})d(0,\partial\Omega)}\,,\, \text{ and }\\ 
&(1-\tilde{t}e^{\tilde{t}})2000|z|\leq |a|\leq (1+\tilde{t}e^{\tilde{t}})2000|z|
\,.
\end{align*}
Especially, it clearly follows that $|z|<|a|<d(0,\partial\Omega)$, $J[z,a]\subset\Omega$, and 
\begin{equation*}
d_{Q}(z,a)\leq \frac{|z-a|}{d(0,\partial\Omega)-|a|}\,.
\end{equation*}

Furthermore, by the elementary estimate from the proof of Proposition \ref{distcurv}, it holds that
\begin{equation*}
|z-a|\leq |a|-\frac{a}{|a|}\cdot z+\frac{|z|^2}{2|a|}\,.
\end{equation*}
Moreover, (\ref{ddoletus}) and (\ref{ekakeissi}) combined with the fact that $\gamma'(t)/|\gamma'(t)|$ is $1$-Lipschitz, imply that
\begin{equation*}
\frac{a}{|a|}\cdot z\geq \frac{1}{100}|z|+|z|\tilde{t}\,\geq \frac{1}{200}|z|\,.
\end{equation*}
 Therefore, we get from above that
\begin{equation*}
|z-a|\leq |a|-\frac{1}{200}|z|+\frac{|z|^2}{2|a|}\leq |a|-\frac{1}{300}|z|\,,
\end{equation*}
where the final inequality follows by the fact $\frac{|z|}{|a|}\approx\frac{1}{2000}$, which is verified above. For the desired contradiction, it suffices to check that 
\begin{equation*}
\frac{|a|-\frac{1}{300}|z|}{d(0,\partial\Omega)-|a|}<
\frac{|a|}{(1+\tilde{t}e^{\tilde{t}})d(0,\partial\Omega)}
\,\,\,\,\,\,\,(
\,\leq\tilde{t}\,).
\end{equation*}
It follows from above that this holds if 
\begin{align*}
&\frac{|a|-\frac{1}{300}|z|}{1-\tilde{t}e^{\tilde{t}}}<\frac{|a|}{1+\tilde{t}e^{\tilde{t}}}\,\\
\Longleftrightarrow\,\,\,\,\,&2|a|\tilde{t}e^{\tilde{t}}<\frac{|z|}{300}(1+\tilde{t}e^{\tilde{t}})\,.
\end{align*}
This follows, since 
\begin{equation*}
\frac{2|a|\tilde{t}e^{\tilde{t}}}{|z|(1+\tilde{t}e^{\tilde{t}})}\leq 2(2000)\tilde{t}e^{\tilde{t}}\leq 10^4\tilde{t}=\frac{10^4(2000)|z|}{d(0,\partial\Omega)}<\frac{1}{300}\,,
\end{equation*}
where the final inequality follows directly from the assumption $|z|\leq 2c(\Omega,x,r,\tilde{r})\,.$

An analogous reasoning implies the remaining case 
\begin{equation*}
\frac{\gamma_{x,y}'(r)}{|\gamma_{x,y}'(r)|}\cdot\frac{z}{|z|}\leq -\frac{1}{100}\,.
\end{equation*} 
In this case one does not need to consider the prolongation of $\gamma$. Indeed, one may choose $a=\gamma(r-\tilde{t})$ and estimate as above that $d_Q(z,a)<\tilde{t}\,$.
\hfill$\Box$\\

In the following corollaries $\Omega, x,r$ and $\tilde{r}$  are as in the previous lemma. 

\begin{corollary}\label{nojoon}
Suppose that $y_1,y_2,y_3\in\partial B_Q(x,r)$ are such that
\begin{equation*}
|y_1-y_2|+|y_2-y_3|\leq 2c(\Omega,x,r,\tilde{r})\,.
\end{equation*}
Then it holds that
\begin{equation}\label{aladdin}
\bigg|\frac{(y_i-y_j)}{|y_j-y_i|}\cdot \frac{\gamma'_{x,y_1}(r)}{|\gamma'_{x,y_1}(r)|}\bigg|\leq\frac{1}{30}\,
\end{equation}
for all $i,j\in\{1,2,3\}$, $i\neq j\,$.
 \end{corollary}
\textit{Proof.} If $i=1$ or $j=1$ in (\ref{aladdin}), then the claim follows directly from Proposition \ref{eiviitetta}. For the remaining case, let us denote $\frac{\gamma_{x,y_1}'(r)}{|\gamma_{x,y_1}'(r)|}=:a$ and $\frac{\gamma_{x,y_2}'(r)}{|\gamma_{x,y_2}'(r)|}=:b$. Since $a\cdot\frac{(y_2-y_1)}{|y_2-y_1|}\leq \frac{1}{100}$ and $b\cdot \frac{(y_2-y_1)}{|y_2-y_1|}\leq\frac{1}{100}$, it follows that
\begin{equation*}
|a-b|\leq \frac{1}{50}\,\text{ or }\,|a+b|\leq \frac{1}{50}\,.
\end{equation*}
This implies the claim, since 
\begin{align*}
\bigg{|}\frac{(y_2-y_3)}{|y_2-y_3|}\cdot a\bigg{|}&=\bigg{|}\mp\frac{(y_2-y_3)}{|y_2-y_3|}\cdot b+(a\pm b)\cdot\frac{(y_2-y_3)}{|y_2-y_3|}\bigg{|}\\
&\leq \bigg{|}\frac{(y_2-y_3)}{|y_2-y_3|}\cdot b\bigg|+\min\{|a+b|,|a-b|\}\,\\
&\leq \frac{1}{100}+\frac{1}{50}\,<\frac{1}{30}\,.
\end{align*}
In the final inequality we again used Proposition \ref{eiviitetta}. 
\hfill$\Box$

\begin{corollary}
Suppose that $\Omega, x,r,$ and $\tilde{r}$ are as before and let 
$y,z\in \partial B_Q(x,r)$ be such that $|z-y|\leq c(\Omega,x,r,\tilde{r})\,.$
Then there exists a $2$-Lipschitz mapping 
$\phi:[0,|z-y|]\to \partial B_Q(x,r)$ such that $\phi(0)=y$ and $\phi(|z-y|)=z\,$.
\end{corollary}
\textit{Proof.}
Using a suitable rotation and translation, we may assume that $y=0$, $z$ is real and $Re(z)>0$. For every $\lambda\in [0,Re(z)]$, let $\phi(\lambda)$ be a point in $\partial B_Q(x,r)\cap B(0,2|z-y|)$ such that $Re(\phi(\lambda))=\lambda$. The existence and the uniqueness of $\phi(\lambda)$ and the desired Lipschitz-property follow easily from Corollary \ref{nojoon} above.



\subsection*{Proof of Lemma \ref{loppusjo3}:} Let $\phi$ be as in the previous proposition and define $G:[0,|z-y|]\to \re$ by
\begin{equation*}
G(\lambda)= \max_{t\in[0,r]}|\gamma_{x,\phi(0)}(t)-\gamma_{x,\phi(\lambda)}(t)\,|\,.
\end{equation*}
For simplicity, let us denote $\gamma_{x,y}=\gamma_y\,$.
It is easy to check that $G$ satisfies 
\begin{align*}
|G(\lambda_1)-G(\lambda_2)|&
=\bigg|\,\max_{t\in[0,r]}|\gamma_{\phi(0)}(t)-\gamma_{\phi(\lambda_1)}(t)\,|-\max_{t\in[0,r]}|\gamma_{\phi(0)}(t)-\gamma_{\phi(\lambda_2}(t)\,|\,\bigg|\\
&\leq \max_{t\in[0,r]}\big|\gamma_{\phi(\lambda_1)}(t)-\gamma_{\phi(\lambda_2)}(t)\,\big|\,.
\end{align*}
For the claim, it suffices to show that $G$ is $2e^{2r}$-Lipschitz, which in turn will follow if for every $\lambda\in[0,|z-y|]$ there exists $\delta_{\lambda}>0$ such that
\begin{equation}\label{foil}
|G(\tilde{\lambda})-G(\lambda)|\leq 2e^{2r}|\tilde{\lambda}-\lambda|\,,\,\text{ if }|\tilde{\lambda}-\lambda|\leq \delta_{\lambda}\,.
\end{equation}
To show (\ref{foil}), recall Corollary \ref{major3}, which implies that for every $y\in\partial B_Q(x,r)$ there exists $\delta_y>0$ such that 
\begin{equation*}
\max_{t\in[0,r]}|\gamma_{\tilde{y}}(t)-\gamma_y(t)|\leq e^{2r}|\tilde{y}-y|\,\text{ if }\tilde{y}\in\partial B_Q(x,r)\text{ and }|\tilde{y}-y|\leq \delta_{y}\,.
\end{equation*}
Combining this with the $2$-Lipschitz property for $\phi$ implies (\ref{foil}). Thus, $G$ is $2e^{2r}$-Lipschitz and the proof is complete. 
\hfill$\Box$\\



 
\subsection*{Proof of Lemma \ref{loppusjo}:}
Suppose that $\Omega$, $r, \tilde{r}$ and $x$ are as in the statement of the lemma. Let again 
$y,z\in \partial B_Q(x,r)$ be such that
 $|y-z|\leq c(\Omega,x,r,\tilde{r})$.
Then it follows from Lemma \ref{loppusjo3} that 
\begin{equation}\label{ekaeka}
\max_{t\in[0,r]}\big|\gamma_{x,y}(t)-\gamma_{x,z}(t)\big|\leq 2e^{2r}|y-z|\,.
\end{equation}
Moreover, it is easy to check that $c(\Omega,x,r,\tilde{r})\leq\frac{m(\Omega,x,r)}{2e^{2r}}$. This yields that Lemma \ref{tili} can be applied with $C=2e^{2r}$ to obtain
\begin{equation*}
d_Q(x,\frac{y+z}{2})\leq r+\frac{4re^{4r}}{(m(\Omega,x,r))^2}|y-z|^2\,\leq r+\frac{4re^{4r}}{(m(\Omega,x,\tilde{r}))^2}|y-z|^2\,.
\end{equation*} 
\hfill$\Box$\\
Then we are ready to prove the main theorems.
\begin{theorem}\label{vihdoin}
Suppose that $\Omega\subset\re^2$ is a domain, $x\in\Omega$ and $r_0<\pi$. Then there exists a unique $\gamma_{x,y}$ for all $y\in B_Q(x,r_0)\,$.
\end{theorem}
\textit{Proof.}
Suppose, on the contrary, that there exists $y_0\in B_Q(x,r_0)$ with two different geodesics $\gamma$ and $\tilde{\gamma}$ from $x$ to $y_0$. By a possible redefinition, we may assume that $d_Q(x,y_0)=r_0$ and (in the same way as before) that $\gamma$ and $\tilde{\gamma}$ are disjoint except the common endpoints. Let us choose a sequence of Voronoi domains $\Omega_k\to\Omega$ as $k\to\infty$ (in the sense expressed in the previous section). Also denote by $B_k(a,t)$ the quasihyperbolic balls with respect to $\Omega_k$ and $B_Q(a,t)$ the quasihyperbolic balls with respect to $\Omega$. Analogously, let $d_k(a,b)$ refer to the quasihyperbolic distance with respect to $\Omega_k$, and $d_Q(a,b)$ with respect to $\Omega$. As before, let
\begin{align*}
m&=\inf_{z\in B_Q(x,\pi)}d(z,\partial\Omega)\,\,\,\text{ and }\\ 
m_k&=\inf_{z\in B_k(x,\pi)}d(z,\partial\Omega_k)\,. 
\end{align*}

Let us then pick $r\in(0,r_0)$ to be close enough to $r_0$ so that 
\begin{equation}\label{rajalle1}
0<|\gamma(r)-\tilde{\gamma}(r)|=:|y-z|\leq \frac{1}{2}C(\Omega,x,r,\pi)\,. 
\end{equation}
Since $\Omega_k\to\Omega$, we get that $m_k\to m$ implying that $C(\Omega_k,x,r,\pi)\to C(\Omega,x,r,\pi)$ as $k\to\infty$. Moreover, we can choose $y_k,z_k\in\partial B_k(x,r)$ such that $y_k\to y$ and $z_k\to z$ as $k\to\infty$. Then it clearly follows that 
\begin{equation}\label{rajalle2}
|y_k-z_k|\leq C(\Omega_k,x,r,\pi)
\end{equation}
when $k$ is large enough. Moreover, by Corollary \ref{Voronoissa1} from the previous section, we know that $\gamma_{x,y}$ (with respect to $\Omega_k$) exists for all $y\in B_k(x,\pi)$. This, together with (\ref{rajalle2}), guarantees that for every $k$ large enough Lemma \ref{loppusjo} can be applied to obtain  
\begin{equation*}
d_k(x,\frac{y_k+z_k}{2})\leq r+\frac{4re^{4r}}{m_k^2}|y_k-z_k|^2\,. 
\end{equation*}
By applying $\Omega_k\to\Omega$, $y_k\to y$ and $z_k\to z$ as $k\to\infty$, we get that
\begin{equation}\label{ref1}
d_Q(x,\frac{y+z}{2})\leq r+\frac{4re^{4r}}{m^2}|y-z|^2\,. 
\end{equation} 

The final contradiction follows again (as in Lemma \ref{versio}) by using Theorem \ref{smallballs}, which implies that
\begin{equation}\label{ref2}
d_{Q}\big(\frac{y+z}{2},y_0\big)
\leq (r_0-r)-\frac{|y-z|^2}{512(r_0-r)(\sup_{w\in B_Q(y_0,r_0-r)}d(w,\partial\Omega))^2}\,.
\end{equation} 
Indeed, by choosing $r_0-r$ small enough, we deduce from (\ref{ref1}) and (\ref{ref2}) above that 
\begin{equation*}
r_0=d_Q(x,y_0)\leq d_{Q}\big(x,\frac{y+z}{2}\big)+d_{Q}\big(\frac{y+z}{2},y_0\big)<r_0\,.
\end{equation*}
This completes the proof.  
\hfill$\Box$\\

\begin{theorem}
Suppose that $\Omega\subset\re^2$ is a simply connected domain. Then there exists a unique $\gamma_{x,y}$ for all $x,y\in \Omega$.
\end{theorem}
\textit{Proof.}
The proof follows very much on the same lines as the proof of the previous Theorem \ref{vihdoin}. Indeed, suppose that there exists $x,y\in \Omega$, $d_Q(x,y)=:r_0$, with two different geodesics $\gamma$ and $\tilde{\gamma}$ from $x$ to $y$. Again, we may assume that $\gamma$ and $\tilde{\gamma}$ are disjoint on $(0,r_0)$ (and canonically parametrized). As before, choose a sequence of Voronoi domains $\Omega_k\to\Omega$ as $k\to\infty$ and let $B_Q(a,t),B_k(a,t),d_Q(a,b)$ and $d_k(a,b)$ be defined as in the above proof.
Furthermore, let
\begin{align*}
m=\inf_{z\in B_Q(x,2r)}d(z,\partial\Omega)\,\,\,\text{ and }\,\,\,\,
m_k=\inf_{z\in B_k(x,2r)}d(z,\partial\Omega_k)\,. 
\end{align*}
By $\Omega_k\to\Omega$ as $k\to\infty$ it again follows that $m_k\to m$ and $C(\Omega_k,x,r_0,2r_0)\to C(\Omega,x,r_0,2r_0)$. Again we pick $r\in(0,r_0)$ so that 
\begin{equation}\label{rajalle12}
0<|\gamma(r)-\tilde{\gamma}(r)|=:|y-z|\leq \frac{1}{2}C(\Omega,x,r_0,2r_0)\,,
\end{equation}
and choose $y_k,z_k\in\Omega_k$ such that $d_k(x,y_k)=d_k(x,z_k)=r$ and $y_k\to y,z_k\to z$ as $k\to\infty$. Moreover, since $\Omega_k\to\Omega$ and $\Omega$ is simply connected, it is possible to apply Corollary \ref{Voronoissa} to obtain the existence of $\gamma_{x,a}$ (with respect to $\Omega_k$) for all $a\in B(x,2r_0)$. Combining this with (\ref{rajalle12}) implies that for every $k$ large enough Lemma \ref{loppusjo} can be applied to obtain  
\begin{equation}\label{again}
d_k(x,\frac{y_k+z_k}{2})\leq r+\frac{4re^{4r}}{m_k^2}|y_k-z_k|^2\,. 
\end{equation}
Then the final contradiction can be obtained by using Theorem \ref{smallballs} exactly in the same way as in the proof of the previous theorem.
\hfill$\Box$\\

\end{document}